\documentclass[12pt]{article}
\usepackage{latexsym,amsfonts,amsthm,amsmath,amscd,amssymb}
\usepackage{graphicx,enumerate,fancyhdr,caption2}

\newtheorem{lem}{Lemma}[section]
\newtheorem{thm}[lem]{Theorem}
\newtheorem{rem}[lem]{Remark}
\newtheorem{prop}[lem]{Proposition}

\newtheorem{conj}[lem]{Conjecture}

\newcommand{\color}[2][1]{}

\newcommand{\dimo}[1]{\vspace{4pt}\noindent\textit{Proof of \ref{#1}}.\ }
\newcommand{\finedimo}{{\hfill\hbox{$\square$}\vspace{4pt}}}
\newcommand{\myprod}{\!\cdot\!}

\newcommand\matT{{\mathbb{T}}}
\newcommand\matN{{\mathbb{N}}}
\newcommand\matZ{{\mathbb{Z}}}
\newcommand\matS{{\mathbb{S}}}

\newcommand\matC{{\mathbb{C}}}
\newcommand\matE{{\mathbb{E}}}
\newcommand\matH{{\mathbb{H}}}
\newcommand\matR{{\mathbb{R}}}
\newcommand\matX{{\mathbb{X}}}

\newcommand\ftil{{\widetilde f}}
\newcommand\ptil{{\widetilde p}}
\newcommand\qtil{{\widetilde q}}

\newcommand\Atil{{\widetilde A}}
\newcommand\Btil{{\widetilde B}}
\newcommand\Ctil{{\widetilde C}}
\newcommand\Dtil{{\widetilde D}}
\newcommand\util{{\widetilde u}}
\newcommand\pitil{{\widetilde \pi}}

\newcommand{\myto}{\mathop{\longrightarrow}\limits}

\renewcommand{\hbar}{{\overline{h}}}

\newfont{\Got}{eufm10 scaled 1200}

\newcommand{\compo}{\,{\scriptstyle\circ}\,}

\newcommand{\mycap} [1] {\caption{\footnotesize{#1}}}

\newcommand{\chiorb}{\chi^{{\mathrm{orb}}}}
\newcommand\Isomp{{\rm Isom}^+}
\newcommand{\Sigmatil}{\widetilde\Sigma}
\newcommand{\Xtil}{\widetilde X}
\newcommand{\Gammatil}{\widetilde\Gamma}

\newcommand{\dotsto}{\mathop{\dashrightarrow}\limits}
\newcommand{\dotstobis}{\mathop{\dashrightarrow\dashrightarrow}\limits}

\newcommand{\argdotsto}[2]{\mathop{\dashrightarrow}\limits^{#1}_{\scriptscriptstyle{#2}}}
\newcommand{\argdotstobis}[2]{\mathop{\dashrightarrow\dashrightarrow}\limits^{#1}_{\scriptscriptstyle{#2}}}
\newcommand{\argdotstoter}[2]{\mathop{\dashrightarrow\dashrightarrow\dashrightarrow}\limits^{#1}_{\scriptscriptstyle{#2}}}
\newcommand{\argdotstoqua}[2]{\mathop{\dashrightarrow\dashrightarrow\dashrightarrow\dashrightarrow}\limits^{#1}_{\scriptscriptstyle{#2}}}
\newcommand{\argdotstoqui}[2]{\mathop{\dashrightarrow\dashrightarrow\dashrightarrow\dashrightarrow\dashrightarrow}\limits^{#1}_{\scriptscriptstyle{#2}}}
\newcommand{\argdotstosex}[2]{\mathop{\dashrightarrow\dashrightarrow\dashrightarrow\dashrightarrow\dashrightarrow\dashrightarrow}\limits^{#1}_{\scriptscriptstyle{#2}}}
\newcommand{\argdotstoset}[2]{\mathop{\dashrightarrow\dashrightarrow\dashrightarrow\dashrightarrow\dashrightarrow\dashrightarrow\dashrightarrow}\limits^{#1}_{\scriptscriptstyle{#2}}}

\def\nequiv{\not\equiv}

\newcommand\calA{{\mathcal A}}

\begin{document}

\title{Surface branched covers\\ and geometric 2-orbifolds}

\author{Maria Antonietta~\textsc{Pascali}
\and\addtocounter{footnote}{5} Carlo~\textsc{Petronio}}

\maketitle

\begin{abstract}
\noindent Let $\Sigmatil$ and $\Sigma$ be closed, connected, and
orientable surfaces, and let $f:\Sigmatil\to\Sigma$ be a branched
cover. For each branching point $x\in\Sigma$ the set of local
degrees of $f$ at $f^{-1}(x)$ is a partition of the total degree
$d$. The total length of the  various partitions is determined by
$\chi(\Sigmatil)$, $\chi(\Sigma)$, $d$ and the number of branching
points via the Riemann-Hurwitz formula. A very old problem asks
whether a collection of partitions of $d$ having the appropriate
total length (that we call a \emph{candidate cover})
always comes from some branched cover. The answer
is known to be in the affirmative whenever $\Sigma$ is not the
$2$-sphere $S$,
while for $\Sigma=S$ exceptions do occur. A
long-standing conjecture however asserts that
when the degree $d$ is a prime number a candidate cover is always realizable.
In this paper we analyze
the question from the point of view of the geometry of
2-orbifolds, and we provide strong supporting evidence for the
conjecture. In particular, we exhibit three different sequences of
candidate covers, indexed by their degree, such that for each sequence:
\begin{itemize}
\item The degrees giving realizable covers have asymptotically
zero density in the naturals;
\item Each prime degree gives a realizable cover.
\end{itemize}
\end{abstract}

\section*{Introduction}
This paper is devoted to a geometric approach, based on
2-orbifolds, to the Hurwitz existence problem for branched covers
between surfaces, see~\cite{Hurwitz} for the original
source, the
classical~\cite{Husemoller,Thom,Singerman,Endler,Francis,Ezell,Mednykh1,EKS,Gersten,Mednykh2,KZ},
the more
recent~\cite{Baranski,MSS,OkPand,PePe1,PePe2,Pako,Zheng}, and below.
We determine the realizability of all candidate surface branched covers
inducing candidate covers between $2$-orbifold with non-negative Euler
characteristic. This yields a complete analysis of the
existence of several infinite series of candidate covers and in
particular to theorems giving rather surprising
connections with number-theoretic facts. These results provide
in particular strong support for the
long-standing conjecture~\cite{EKS} that a candidate cover with
prime total degree is always realizable. The appropriate
terminology and the statements of some of the results
established below are given in
the rest of the present Introduction. Our main results
(Theorems~\ref{244:main:thm} to~\ref{333:main:thm}) will not
involve explicit reference to 2-orbifolds.

\paragraph{Branched covers}
Let $\Sigmatil$ and $\Sigma$ be closed, connected, and orientable surfaces, and
$f:\Sigmatil\to\Sigma$ be a branched cover, \emph{i.e.} a map
locally modelled on functions of the form $(\matC,0)\myto^{z\mapsto z^k}(\matC,0)$
with $k\geqslant 1$. If $k>1$ then $0$ in the target $\matC$ is  a
branching point, and $k$ is the local degree at $0$ in the
source $\matC$. There is a finite number $n$ of branching points,
and, removing all of them from $\Sigma$ and their preimages
from $\Sigmatil$, we see that $f$ induces a genuine cover
of some degree $d$. The collection $(d_{ij})_{j=1}^{m_i}$ of
the local degrees at the preimages of the $i$-th branching point is a
partition $\Pi_i$ of $d$.
We now define:
\begin{itemize}
\item $\ell(\Pi_i)$ to be the length $m_i$ of $\Pi_i$;
\item $\Pi$ as the set $\{\Pi_1,\ldots,\Pi_n\}$ of all partitions of $d$ associated to $f$;
\item $\ell(\Pi)$ to be the total length $\ell(\Pi_1)+\ldots+\ell(\Pi_n)$ of $\Pi$.
\end{itemize}
Then multiplicativity of the Euler characteristic $\chi$ under
genuine covers for surfaces with boundary implies
the classical Riemann-Hurwitz formula
\begin{equation}\label{RH:general:eq}
\chi(\Sigmatil)-\ell(\Pi)=d\myprod\big(\chi(\Sigma)-n\big).
\end{equation}

\paragraph{Candidate branched covers and the realizability problem}
Consider again two closed, connected, and orientable surfaces $\Sigmatil$ and $\Sigma$,
integers $d\geqslant 2$ and $n\geqslant 1$, and a set of partitions
$\Pi=\{\Pi_1,\ldots,\Pi_n\}$ of $d$, with $\Pi_i=(d_{ij})_{j=1}^{m_i}$, such
that condition~(\ref{RH:general:eq}) is satisfied. We associate to these data
the symbol
$$\Sigmatil\argdotstosex{d:1}{(d_{11},\ldots,d_{1m_1}),\ldots,(d_{n1},\ldots,d_{nm_n})}\Sigma$$
that we will call a \emph{candidate surface branched cover}. A classical (and
still not completely solved) problem, known as the \emph{Hurwitz existence problem},
asks which candidate surface branched covers are actually \emph{realizable},
namely induced by some existent branched cover $f:\Sigmatil\to\Sigma$.
A non-realizable candidate surface branched cover will be called \emph{exceptional}.

\begin{rem}
\emph{The symbol $\big(\Sigmatil,\Sigma,n,d,(d_{ij})\big)$ and the name
\emph{compatible branch datum} are used in \cite{PePe1,PePe2} instead
of the terminology of ``candidate covers'' we will use here.}
\end{rem}

\paragraph{Known results}
Over the last 50 years the Hurwitz existence problem was the object of many papers,
already  listed above. The combined efforts of several
mathematicians led in particular to the following
results~\cite{Husemoller,EKS}:
\begin{itemize}
\item If $\chi(\Sigma)\leqslant 0$ then any candidate surface branched cover is realizable,
\emph{i.e.} the Hurwitz existence problem has a positive solution in this case;
\item If $\chi(\Sigma)>0$, \emph{i.e.} if $\Sigma$ is the $2$-sphere $S$, there exist
exceptional candidate surface branched covers.
\end{itemize}

\begin{rem}
\emph{A version of the Hurwitz existence problem exists also for
possibly non-orientable $\Sigmatil$ and $\Sigma$. Condition~(\ref{RH:general:eq})
must be complemented in this case with a few more requirements
(some obvious, and one slightly less obvious, see~\cite{PePe1}).
However it has been shown~\cite{Ezell,EKS}
that again this generalized problem always has a positive solution if $\chi(\Sigma)\leqslant 0$,
and that the case where $\Sigma$ is the projective plane reduces to the case
where $\Sigma$ is the $2$-sphere $S$.}
\end{rem}

According to these facts, in order to face the Hurwitz existence problem,
it is not restrictive to \emph{assume the candidate covered surface $\Sigma$ is the $2$-sphere $S$},
which we will do henceforth.
Considerable energy has been devoted over the time to a general understanding of
the exceptional candidate surface branched covers in this case, and
quite some progress has been made (see for instance the survey of known results
contained in~\cite{PePe1}, together with the later papers~\cite{PePe2,Pako,Zheng}),
but the global pattern remains elusive. In particular the following conjecture
proposed in~\cite{EKS} appears to be still open:

\begin{conj}\label{prime:conj}
If $\Sigmatil\argdotsto{d:1}{\Pi}S$ is candidate surface branched cover and
the degree $d$ is a prime number then the candidate is realizable.
\end{conj}

We mention in particular that all exceptional candidate surface branched covers
with $n=3$ and $d\leqslant 20$
have been determined by computer in~\cite{Zheng}. There are very many of them, but
none occurs for prime $d$.

\paragraph{Main new results}
Using the geometry of Euclidean 2-orbifolds we will establish in this paper
(among others) the next three theorems. Recall that a candidate surface branched cover
is a set of data $\Sigmatil\argdotstobis{d:1}{\Pi_1,\ldots,\Pi_n}\Sigma$
satisfying the Riemann-Hurwitz condition~(\ref{RH:general:eq}). Moreover $S$ denotes
the 2-sphere.

\begin{thm}\label{244:main:thm}
Suppose $d=4k+1$ for $k\in\matN$. Then
$$S\argdotstosex{d:1}
{(\underbrace{{\scriptscriptstyle 2,\ldots,2}}_{2k},1),
(\underbrace{{\scriptscriptstyle 4,\ldots,4}}_{k},1),
(\underbrace{{\scriptscriptstyle 4,\ldots,4}}_{k},1)}S$$
is a candidate surface branched cover, and
it is realizable if and only if $d$ can be expressed as $x^2+y^2$ for some $x,y\in\matN$.
\end{thm}

\begin{thm}\label{236:main:thm}
Suppose $d=6k+1$ for $k\in\matN$. Then
$$S\argdotstosex{d:1}
{(\underbrace{{\scriptscriptstyle 2,\ldots,2}}_{3k},1),
(\underbrace{{\scriptscriptstyle 3,\ldots,3}}_{2k},1),
(\underbrace{{\scriptscriptstyle 6,\ldots,6}}_{k},1)}S$$
is a candidate surface branched cover and
it is realizable if and only if $d$ can be expressed as $x^2+xy+y^2$ for some $x,y\in\matN$.
\end{thm}

\begin{thm}\label{333:main:thm}
Suppose $d=3k+1$ for $k\in\matN$. Then
$$S\argdotstosex{d:1}
{(\underbrace{{\scriptscriptstyle 3,\ldots,3}}_{k},1),
(\underbrace{{\scriptscriptstyle 3,\ldots,3}}_{k},1),
(\underbrace{{\scriptscriptstyle 3,\ldots,3}}_{k},1)}S$$
is a candidate surface branched cover and
it is realizable if and only if $d$ can be expressed as $x^2+xy+y^2$ for some $x,y\in\matN$.
\end{thm}

What makes these results remarkable in
view of Conjecture~\ref{prime:conj} is that:
\begin{itemize}
\item A prime number of the form $4k+1$ can always be expressed as
$x^2+y^2$ for $x,y\in\matN$ (Fermat);
\item A prime number of the form $6k+1$ (or equivalently $3k+1$)
can always be expressed as
$x^2+xy+y^2$ for $x,y\in\matN$ (Gauss);
\item The integers that can be expressed as $x^2+y^2$ or as
$x^2+xy+y^2$ with $x,y\in\matN$ have asymptotically zero density in $\matN$.
\end{itemize}
This means that a candidate cover in any of our three statements
is ``exceptional with probability 1,'' even though it is
realizable when its degree is prime. Note also that it was
shown in~\cite{EKS} that establishing Conjecture~\ref{prime:conj} in the
special case of three branching points would imply the general case.

\paragraph{Induced candidate 2-orbifold covers}
A 2-orbifold $X=\Sigma(p_1,\ldots,p_n)$ is
a closed orientable surface $\Sigma$
with $n$ cone points of orders $p_i\geqslant 2$, at which $X$ has a singular
differentiable structure given by the quotient
$\matC/_{\langle{\rm rot}(2\pi/p_i)\rangle}$.
Bill Thurston~\cite{thurston:notes}
introduced the notions (reviewed below) of orbifold cover and orbifold Euler
characteristic
$$\chiorb\big(\Sigma(p_1,\ldots,p_n)\big)=\chi(\Sigma)-\sum_{i=1}^n\left(1-\frac1{p_i}\right),$$
designed so that if $f:\Xtil\myto^{d:1} X$ is an orbifold cover then
$\chiorb(\Xtil)=d\cdot\chiorb(X)$. He also showed that:
\begin{itemize}
\item If $\chiorb(X)>0$ then $X$ is either \emph{bad} (not covered by a surface in the sense of
orbifolds) or
\emph{spherical}, namely the quotient of the metric 2-sphere $\matS$ under a finite isometric action;
\item If $\chiorb(X)=0$ (respectively, $\chiorb(X)<0$) then $X$ is
\emph{Euclidean} (respectively, \emph{hyperbolic}), namely
the quotient of the Euclidean plane $\matE$ (respectively, the hyperbolic plane $\matH$)
under a discrete isometric action.
\end{itemize}

As pointed out in~\cite{PePe1} and spelled out below, any
candidate surface branched cover $\Sigmatil\dotsto^{d:1}_\Pi\Sigma$ induces
a \emph{candidate $2$-orbifold cover} $\Xtil\dotsto^{d:1}X$ satisfying the
condition $\chiorb(\Xtil)=d\myprod\chiorb(X)$. Moreover
$\Sigmatil\dotsto^{d:1}_\Pi\Sigma$ can be reconstructed from
$\Xtil\dotsto^{d:1}X$ if some additional \emph{covering instructions} are provided.

\paragraph{More new results}
The main idea of this paper is to analyze the realizability of a given
candidate surface branched cover $\Sigmatil\dotsto^{d:1}_\Pi\Sigma$ using the
associated candidate 2-orbifold cover $\Xtil\dotsto^{d:1}X$
and the geometries of $\Xtil$ and $X$. It turns out that this is particularly
effective when $\chiorb(X)$ is non-negative (note that $\chiorb(\Xtil)$ has the
same sign as $\chiorb(X)$, being $d$ times it), but we will
also briefly touch the case $\chiorb(X)<0$.
In the bad/spherical case the statement is quite expressive:

\begin{thm}\label{pos:chi:main:thm}
Let a candidate surface branched cover $\Sigmatil\dotsto^{d:1}_\Pi\Sigma$
induce a candidate $2$-orbifold cover $\Xtil\dotsto^{d:1}X$ with $\chiorb(X)>0$.
Then $\Sigmatil\dotsto^{d:1}_\Pi\Sigma$ is exceptional if and only if
$\Xtil$ is bad and $X$ is spherical. All exceptions occur with
non-prime degree.
\end{thm}

Turning to the Euclidean case (which leads in particular to
Theorems~\ref{244:main:thm} to~\ref{333:main:thm}) we confine ourselves here
to the following informal:

\begin{thm}\label{zero:chi:main:thm}
Let a candidate surface branched cover $\Sigmatil\dotsto^{d:1}_\Pi\Sigma$
induce a candidate $2$-orbifold cover $\Xtil\dotsto^{d:1}X$ with Euclidean $X$.
Then its realizability can be decided explicitly in terms of $d$ and $\Pi$.
More precisely, as in Theorems~\ref{244:main:thm} to~\ref{333:main:thm},
given $\Pi$ the condition on
$d$ depends on a congruence and/or an integral
quadratic form. No exceptions occur when $d$ is a prime number.
\end{thm}

We conclude with our statement for the hyperbolic case.
A $2$-orbifold is called \emph{triangular} if it has the form $S(p,q,r)$.

\begin{thm}\label{neg:chi:main:thm}
There exist $9$ candidate surface branched covers
inducing a candidate $2$-orbifold cover $\Xtil\dotsto X$
with $\Xtil$ and $X$ being hyperbolic triangular orbifolds.
All of them but two are realizable. Exceptions occur in degrees 8 and $16$ (which
are not prime).
\end{thm}

For the case of non-negative $\chiorb$,
within the proofs of Theorems~\ref{pos:chi:main:thm}
and~\ref{zero:chi:main:thm} we will describe explicit geometric
constructions of all the realizable covers. The
proof of Theorem~\ref{neg:chi:main:thm} has instead a more combinatorial
flavour.

\vspace{.5cm}

\noindent\textsc{Acknowledgements}. Part of this work was carried
out while the second named author was visiting the Universit\'e
Paul Sabatier in Toulose and the Columbia University in New York.
He is grateful to both these institutions for financial support,
and he would like to thank Michel Boileau and Dylan Thurston for
their warm hospitality and inspiring mathematical discussions.

\section{A geometric approach to the problem}\label{approach:sec}
In this section we will describe the general framework leading to the
proofs of Theorems~\ref{244:main:thm} to~\ref{zero:chi:main:thm}, carried out
in Sections~\ref{pos:chi:sec} and~\ref{zero:chi:sec}.

\paragraph{Orbifold covers}
Besides the terminology and facts on $2$-orbifolds already mentioned
in the Introduction (and not reviewed here) we will need the precise
definition of a degree-$d$ cover $f:\Xtil\to X$ between 2-orbifolds.
This is a map
such that $f^{-1}(x)$ generically consists of $d$ points and
locally making a diagram of the following form commutative:
$$\begin{array}{ccc}
(\matC,0) & \myto^{\rm id} & (\matC,0) \\
\downarrow & & \downarrow \\
(\Xtil,\widetilde{x}) & \myto^{f} & (X,x)
\end{array}$$
where $\widetilde{x}$ and $x$ have cone orders $\widetilde{p}$
and $p=k\myprod\widetilde{p}$ respectively, and the vertical arrows are the
projections corresponding to the actions of
$\langle{\rm rot}(2\pi/\widetilde{p})\rangle$ and $\langle{\rm rot}(2\pi/p)\rangle$,
namely the maps defining the (possibly singular) local differentiable structures at
$\widetilde{x}$ and $x$. Since this local model can be described by the map $z\mapsto z^k$,
we see that $f$ induces a branched cover between the underlying surfaces of $\Xtil$ and $X$.
Using the Riemann-Hurwitz formula~(\ref{RH:general:eq})
it is then easy to show that $\chiorb(\Xtil)=d\myprod\chiorb(X)$.

Bill Thurston introduced in~\cite{thurston:notes} the notion of \emph{orbifold universal cover}
and established its existence. Using this fact we easily get the following result that we
will need below (where, not surprisingly, \emph{good} means ``not bad''):

\begin{lem}\label{no:bad:on:good:lem}
If $\Xtil$ is bad and $X$ is good then there cannot exist any orbifold cover $\Xtil\to X$.
\end{lem}

\paragraph{Induced orbifold covers}
As one easily sees, distinct orbifold covers can induce the same surface
branched cover (in the local model, the two cone orders can be multiplied by one and the
same integer). However a surface branched cover has an ``easiest'' associated
orbifold cover, \emph{i.e.} that with the smallest possible cone orders. This carries over
to \emph{candidate} covers, as we will now spell out.
Consider a candidate surface branched cover
$$\Sigmatil\argdotstosex{d:1}{(d_{11},\ldots,d_{1m_1}),\ldots,(d_{n1},\ldots,d_{nm_n})}\Sigma$$
and define
$$\begin{array}{ll}
p_i={\rm l.c.m.}\{d_{ij}:\ j=1,\ldots,m_i\},\quad & p_{ij}=p_i/d_{ij},\phantom{\Big|} \\
X=\Sigma(p_1,\ldots,p_n), & \Xtil=\Sigmatil\big((p_{ij})_{i=1,\ldots,n}^{j=1,\ldots,m_i}\big)
\end{array}$$
where ``l.c.m.'' stands for ``least common multiple.''
Then we have an induced candidate $2$-orbifold cover $\Xtil\dotsto^{d:1}X$
satisfying $\chiorb(\Xtil)=d\myprod\chiorb(X)$. Note that the original
candidate surface branched cover cannot be reconstructed from $\Xtil,X,d$ alone,
but it can if $\Xtil\dotsto^{d:1}X$ is complemented with the \emph{covering instructions}
$$(p_{11},\ldots,p_{1m_1})\dotsto p_1,\qquad\ldots\qquad
(p_{n1},\ldots,p_{nm_n})\dotsto p_n$$
that we will sometimes include in the symbol $\Xtil\dotsto^{d:1}X$ itself,
omitting the $p_{ij}$'s equal to $1$. Of course
a candidate surface branched cover is realizable if and only if the induced
candidate 2-orbifold cover with appropriate covering instructions is realizable.

\paragraph{The geometric approach}
To analyze the realizability of a candidate surface branched cover we will
switch to the induced candidate 2-orbifold cover
$\Xtil\dotsto X$ and we will use geometry either to explicitly
construct a map $f:\Xtil\to X$ realizing it, or to
show that such an $f$ cannot exist.

To explain how this works
we first note that any 2-orbifold $X$ with a fixed
geometric structure of type $\matX\in\{\matS,\matE,\matH\}$
has a well-defined distance function. This is because the structure is given
by a quotient map
$\matX\to X$, that for obvious reasons we will call \emph{geometric
universal cover of $X$}, defined by an isometric and
discrete (even if not free) action. Therefore a piecewise smooth path $\alpha$
in $X$ has a well-defined length obtained by lifting it to a path
$\widetilde{\alpha}$ in $\matX$, even if $\widetilde{\alpha}$ itself is not
unique (even up to automorphisms of $\Xtil$)
when $\alpha$ goes through some cone point of $X$.
Now we have the following:

\begin{prop}\label{geom:X:then:geom:Xtil:prop}
Let $f:\Xtil\rightarrow X$ be a $2$-orbifold cover. Suppose that $X$
has a fixed geometry with geometric
universal cover $\pi:\matX\rightarrow X$. Then there exists a geometric structure on
$\Xtil$ with geometric universal cover
$\pitil:\matX\rightarrow \Xtil$ and an isometry
$\ftil:\matX\rightarrow \matX$ such that
$\pi\compo\ftil=f\compo\pitil$.
\end{prop}

\begin{proof}
We define the length of a path in $\Xtil$ as the length of its image in $X$ under $f$,
and we consider the corresponding distance. Analyzing the local model of $f$,
one sees that this distance is compatible with a local orbifold geometric
structure also of type $\matX$, so there is one global such structure
on $\Xtil$, with geometric universal cover
$\pitil:\matX\rightarrow \Xtil$. The properties of the
universal cover imply that there is a map $\ftil:\matX\rightarrow \matX$ such that
$\pi\compo\ftil=f\compo\pitil$. By construction $\ftil$ preserves
the length of paths, but $\matX$ is a manifold, not an orbifold,
so $\ftil$ is a local isometry. In particular it is a cover, but $\matX$ is simply
connected, so $\ftil$ is a homeomorphism and hence an isometry.
\end{proof}

Any spherical 2-orbifold $X$ is \emph{rigid}, namely the
geometric universal cover $\matS\to X$ is unique up to automorphisms of $\matS$ and
$X$, so in the spherical case one is not faced with any choice while applying
Proposition~\ref{geom:X:then:geom:Xtil:prop}.
On the contrary Euclidean 2-orbifolds are never rigid, since the metric can always be rescaled
(and it can also be changed in more essential ways on the torus $T$ and on $S(2,2,2,2)$, see below).
In this case we will slightly modify the content of
Proposition~\ref{geom:X:then:geom:Xtil:prop} by rescaling $\Xtil$ so that its area equals that of $X$,
in which case $\ftil$ is no more an isometry but merely a complex-affine map $\matC\to\matC$,
with $\matC$ identified to $\matE$. More precisely:

\begin{prop}\label{Eucl:X:then:Eucl:Xtil:prop}
Let $f:\Xtil\rightarrow X$ be a $2$-orbifold cover. Suppose that $X$
has a fixed Euclidean structure with geometric
universal cover $\pi:\matE\rightarrow X$.
Then there exists a Euclidean structure on
$\Xtil$ with geometric universal cover
$\pitil:\matE\rightarrow \Xtil$
such that $X$ and $\Xtil$ have the same area,
and a map
$\ftil:\matE\rightarrow \matE$ of the form $\ftil(z)=\lambda\myprod z+\mu$ such that
$\pi\compo\ftil=f\compo\pitil$ and $d=|\lambda|^2$.
\end{prop}

\begin{proof}
With respect to the structure on $\Xtil$ given by
Proposition~\ref{geom:X:then:geom:Xtil:prop} the area
of $\Xtil$ is $d$ times that of $X$, so the scaling factor is $1/\sqrt{d}$.
After rescaling $\ftil$ is therefore $\sqrt{d}$ times an isometry, and
the conclusion follows.
\end{proof}

To conclude the section, we note that for candidate
covers of the form $S\argdotsto{d:1}{\Pi}S$, as most of ours will be, the
Riemann-Hurwitz formula~(\ref{RH:general:eq}) reads
\begin{equation}\label{RH:SS3:eq}
\ell(\Pi)=d+2.
\end{equation}

\section{Positive Euler characteristic}\label{pos:chi:sec}

In this section we will establish Theorem~\ref{pos:chi:main:thm}. More precisely
we will show:

\begin{thm}\label{bad:on:good:list:thm}
A candidate surface branched cover $\Sigmatil\dotsto\Sigma$ inducing a candidate $2$-orbifold cover
$\Xtil\dotsto X$ with $\chiorb(X)>0$
is exceptional if and only if $\Xtil$ is
bad and $X$ is spherical. This occurs precisely for the following candidate covers,
in none of which the degree is prime:
\begin{equation}\label{bad:on:good:list:eq}
\begin{array}{lll}
S\argdotstoqui{9:1}{(2,\ldots,2,1),(3,3,3),(3,3,3)}S &
S\argdotstoqui{9:1}{(2,\ldots,2,1),(3,3,3),(4,4,1)}S &
S\argdotstoqui{10:1}{(2,\dots,2),(3,3,3,1),(4,4,2)}S \\
S\argdotstoqui{16:1}{(2,\dots,2),(3,\dots,3,1),(4,\ldots,4)}S &
S\argdotstoqui{16:1}{(2,\dots,2),(3,\dots,3,1),(5,5,5,1)}S &
S\argdotstoqui{18:1}{(2,\dots,2),(3,\dots,3),(4,\ldots,4,2)}S \\
S\argdotstoqui{21:1}{(2,\dots,2,1),(3,\dots,3),(5,\ldots,5,1)}S &
S\argdotstoqui{25:1}{(2,\dots,2,1),(3,\dots,3,1),(5,\ldots,5)}S &
S\argdotstoqui{36:1}{(2,\dots,2),(3,\dots,3),(5,\dots,5,1)}S\\
S\argdotstoqui{40:1}{(2,\dots,2),(3,\dots,3,1),(5,\dots,5)}S &
S\argdotstoqui{45:1}{(2,\dots,2,1),(3,\dots,3),(5,\dots,5)}S &
S\argdotstoqui{2k:1}{(2,\dots,2),(2,\dots,2),(h,2k-h)}S
\end{array}
\end{equation}
with $k>h\geqslant 1$ in the last item.
\end{thm}

In addition to proving this result we will describe
all $\Sigmatil\dotsto\Sigma$
inducing $\Xtil\dotsto X$ with $\chiorb(X)>0$
not listed in the statement, and we will
explicitly construct a geometric realization of each such $\Xtil\dotsto X$.
To outline our argument, we first recall
that the $2$-orbifolds $X$ with $\chiorb(X)>0$ are
$$S,\quad  S(p),\quad S(p,q),\quad S(2,2,p),\quad S(2,3,3),
\quad S(2,3,4),\quad S(2,3,5).$$
In particular for any relevant $\Sigmatil\argdotsto{d:1}{\Pi}\Sigma$
we have $\Sigmatil=\Sigma=S$.
Moreover $X$ is bad if and only if it is $S(p)$ for $p>1$ or $S(p,q)$ for
$p\neq q>1$, and in all other cases it has a rigid spherical structure.
Our main steps will be as follows:

\begin{itemize}

\item We will determine all
candidate surface branched covers inducing some $\Xtil\dotsto X$
with positive $\chiorb(\Xtil)$ and $\chiorb(X)$, and the corresponding
covering instructions for $\Xtil\dotsto X$;

\item For each spherical $X$ with $\chiorb(X)>0$ we will explicitly describe
(and fix) the geometric universal cover $\pi:\matS\to X$;

\item For each $\Xtil\dotsto X$ with $\chiorb(X)>0$
(complemented with its covering instructions)
induced by some candidate surface branched cover,
except when $\Xtil$ is bad and $X$ is spherical,
we will explicitly describe an isometry $\ftil:\matS\to\matS$
such that there exists $f:\Xtil\to X$ realizing $\Xtil\dotsto X$
with $\pi\compo\ftil=f\compo\pitil$, where $\pi$ and $\pitil$
are the geometric universal covers of $\Xtil$ and $X$ described in the previous step.
\end{itemize}

\paragraph{Relevant candidate covers}
To list all candidate surface branched covers inducing
some $\Xtil\dotsto X$ with positive $\chiorb(X)$
our steps will be as follows:
\begin{itemize}
\item We consider all possible pairs
($\Xtil,X)$ such that $\chiorb(\Xtil)/\chiorb(X)$ is an integer $d>1$;
\item Supposing $X$ has $n$ cone points of orders $p_1,\ldots,p_n$, we consider
all possible ways of grouping the orders of the cone points of $\Xtil$ as
$$(q_{11},\ldots,q_{1\mu_1}),\ldots,(q_{n1},\ldots,q_{n\mu_n})$$
so that $q_{ij}$ divides $p_i$ for all $i$ and $j$;
\item We determine $m_i\geqslant \mu_i$ so that, setting $q_{ij}=1$ for $j>\mu_i$
and $d_{ij}=\frac{p_i}{q_{ij}}$, we have that $\sum\limits_{j=1}^{m_i}d_{ij}$ is
equal to $d$ for all $i$;
\item We check that $p_i$ is the least common multiple of $(d_{ij})_{j=1}^{m_i}$.
\end{itemize}
This leads to the candidate surface branched cover
$S\argdotsto{d:1}{\Pi}S$ with
$\Pi_i=(d_{ij})_{j=1}^{m_i}$ and $\Pi=(\Pi_i)_{i=1}^n$,
inducing $\Xtil\dotsto X$ with covering instructions
$$(q_{11},\ldots,q_{1m_1})\dotsto p_1,\ldots,(q_{n1},\ldots,q_{nm_n})\dotsto p_n.$$
For the sake of brevity we will group together our statements depending on the type of
$\Xtil$. In the proofs it will sometimes be convenient to carry out the steps outlined
above in a different order. In particular, it is often not easy to determine beforehand when
$\chiorb(\Xtil)/\chiorb(X)$ is an integer, so this condition is imposed at the end. Moreover,
whenever $X$ has three cone points, instead of
$\chiorb(\Xtil)=d\myprod\chiorb(X)$ we will use the equivalent
formula~(\ref{RH:SS3:eq}), expressed in terms of the data of the would-be candidate
surface branched cover.

\begin{rem}\label{triv:cand:rem}
\emph{Suppose we fix some $\Xtil\dotsto^{d:1} X$ where $X$ has $n\geqslant 1$ cone points
and some $\Sigmatil\dotstobis^{d:1}_{\Pi_1,\ldots,\Pi_n}\Sigma$
induces it. Then each of the $n$ partitions of $d$ in $\Pi$
has at least one entry larger than $1$, otherwise $\Sigmatil\dotsto^{d:1}_{\Pi}\Sigma$
induces some $\widetilde{Y}\dotsto Y$ where $Y$ has less than $n$ cone points.
In particular $d>1$.}
\end{rem}

\begin{prop}\label{bad:Xtil:cand:list:prop}
The candidate surface branched covers inducing some
$\Xtil\dotsto X$ with bad $\Xtil$ are precisely those listed
in~(\ref{bad:on:good:list:eq}).
\end{prop}

\begin{proof}
We start with $\Xtil=S(\ptil)$ for $\ptil\geqslant 2$. If $X=S$ or $X=S(p)$ then
$\chiorb(\Xtil)/\chiorb(X)<2$, so there is no relevant candidate.

Now suppose $X=S(p,q)$ and $\ptil\dotsto p$, so $p=k\myprod\ptil$ for some
$k$, whence $\mu_1=1$, $\mu_2=0$, and $d=k+(m_1-1)\myprod p=m_2\myprod q$.
Combining these relations with $1+\frac{1}{\ptil}=d\myprod\left(\frac1p+\frac1q\right)$
we get $m_1+m_2=2$, so $m_1=m_2=1$ and $k=q=d$, but $p$ is not l.c.m.$(k)$, so again
there is no relevant candidate.

Turning to $X=S(2,2,p)$ we can have either $\ptil\dotsto 2$ or $\ptil\dotsto p$.
In the first case we should have $\Pi_1=(2,\ldots,2,1)$ and $\Pi_2=(2,\ldots,2)$, which is
impossible because $d$ should be both even and odd. In the second case we have $d=2k$ and $m_1=m_2=k$, whence $m_3=2$,
so we get item 12 in~(\ref{bad:on:good:list:eq})
in the special case where $h$ divides $2k-h$ or conversely.

Now let $X=S(2,3,3)$. If $\ptil\dotsto 2$ (or $\ptil\dotsto 3$) then $\ptil=2$ (or $\ptil=3$)
and computing $\chiorb$ we get $d=9$ (or $d=8$). In the first case we get
item 1 in~(\ref{bad:on:good:list:eq}), in the second case we get nothing because
$\Pi_2=(3,\ldots,3,1)$ which is incompatible with $d=8$.

The discussion for $X=S(2,3,p)$ with $p=4,5$ is similar. We examine where $\ptil$
can be mapped to, we deduce what it is
(except that both $2$ and $4$ are possible when $\ptil\dotsto p=4$), in each case we determine $d$ using $\chiorb$
and we check that there exist appropriate partitions of $d$.
For $X=S(2,3,4)$ we get
items 4 and 6 in~(\ref{bad:on:good:list:eq}), with $\ptil=2\dotsto 4$ in 6, while for
$X=S(2,3,5)$ we get items 9 to 11 in~(\ref{bad:on:good:list:eq}).

Let us now consider $\Xtil=S(\ptil,\qtil)$ with $\ptil\neq\qtil>1$ and again examine the
various $X$'s, noting first that $X$ cannot be $S$ or $S(p)$ since $d\geqslant 2$.
For $X=S(p,q)$ suppose first $\ptil,\qtil\dotsto p$. Then $p=k\myprod\ptil=h\myprod\qtil$
and $d=k+h+(m_1-2)\myprod p=m_2\myprod q$, which we can combine with $\frac1\ptil+\frac1\qtil=d\myprod\left(\frac1p+\frac1q\right)$
easily getting $m_1=2$ and $m_2=0$, which is impossible. Now suppose $\ptil\dotsto p$ and
$\qtil\dotsto q$, so $p=k\myprod\ptil$ and $q=h\myprod\qtil$ whence $d=k+(m_1-1)\myprod p=h+(m_2-1)\myprod q$,
which leads to $m_1=m_2=1$, but then we cannot have $p={\rm l.c.m.}(k)$ or $q={\rm l.c.m.}(h)$, so we get nothing.

If $X=S(2,2,p)$ then we cannot have $(\ptil,\qtil)\dotsto 2$ or
$\ptil\dotsto 2,\qtil\dotsto 2$,
otherwise $X$ would be good.
If $\ptil\dotsto 2$ and $\qtil\dotsto p$ then $d$ should be both even and odd, which is impossible.
So $(\ptil,\qtil)\dotsto p$, $d=2k$, $m_1=m_2=k$, $m_3=2$ and we get the last item in~(\ref{bad:on:good:list:eq})
with $h$ and $2k-h$ not multiple of each other.

For $X=S(2,3,p)$, considering where $\ptil$ and $\qtil$ can be mapped, we again see what
they can be,
we determine $d$ using $\chiorb$, and we check that three appropriate partitions exist, getting nothing
for $p=3$, items 2 and 3 in~(\ref{bad:on:good:list:eq}) for $p=4$,
and items 5, 7, and 8 for $p=5$, which completes the proof.
\end{proof}

\begin{prop}\label{Xtil:Spp:cand:list:prop}
The candidate surface branched covers inducing some
$S(\ptil,\ptil)\dotsto X$ with $\ptil>1$ are
\begin{equation}\label{Xtil:Spp:cand:list:eq}
\begin{array}{lll}
S\argdotstoqui{4:1}{(2,2),(3,1),(3,1)}S &
S\argdotstoqui{6:1}{(2,2,1,1),(3,3),(3,3)}S &
S\argdotstoqui{6:1}{(2,2,2),(3,3),(4,1,1)}S \\
S\argdotstoqui{8:1}{(2,\ldots,2),(3,3,1,1),(4,4)}S &
S\argdotstoqui{12:1}{(2,\dots,2,1,1),(3,3,3,3),(4,4,4)}S &
S\argdotstoqui{12:1}{(2,\dots,2),(3,\ldots,3),(4,4,2,2)}S \\
S\argdotstoqui{12:1}{(2,\dots,2),(3,\ldots,3),(5,5,1,1)}S &
S\argdotstoqui{20:1}{(2,\dots,2),(3,\dots,3,1,1),(5,\dots,5)}S &
S\argdotstoqui{30:1}{(2,\dots,2,1,1),(3,\dots,3),(5,\dots,5)}S \\
S\argdotstoqui{2k+1:1}{(2,\dots,2,1),(2,\dots,2,1),(2k+1)}S &
S\argdotstoqui{2k+2:1}{(2,\dots,2,1,1),(2,\dots,2),(2k+2)}S &
\end{array}
\end{equation}
with arbitrary $k\geqslant 1$ in the last two items.
\end{prop}

\begin{proof}
Since $\chiorb(S(\ptil,\ptil))=\frac2\ptil\leqslant1$ and $d\geqslant 2$ we
cannot have $X=S$ or $X=S(p)$. Suppose then $X=S(p,q)$, so $\frac2\ptil=\frac dp+\frac dq$.
If $\ptil,\ptil\dotsto p$ then $p=k\myprod\ptil$ and $d=2k+(m_1-2)\myprod p=m_2\myprod q$, whence
$m_1=2$ and $m_2=0$, which is absurd. If $\ptil\dotsto p$ and $\ptil\dotsto q$
then $p=k\myprod\ptil$ and $q=h\myprod\ptil$ whence $d=k+(m_1-1)\myprod p=h+(m_2-1)\myprod q$
which gives $m_1=m_2=1$, but then we cannot have $p={\rm l.c.m.}(k)$ or
$q={\rm l.c.m.}(h)$, so we get nothing.

Assume now $X=S(2,2,p)$.
If $(\ptil,\ptil)\dotsto 2$ then $\ptil=2$ and $p=d=2k+2$,
which leads to the last item in~(\ref{Xtil:Spp:cand:list:eq}).
If $\ptil\dotsto 2$ and $\ptil\dotsto 2$ then $\ptil=2$ and $p=d=2k+1$,
so we get the penultimate item in~(\ref{Xtil:Spp:cand:list:eq}).
Of course we cannot have
$\ptil\dotsto 2$ and $\ptil\dotsto p$ otherwise $d$ should be both even and odd.
If $(\ptil,\ptil)\dotsto p$ then
$p=k\myprod\ptil$ and $d=2k$, so $\Pi_3=(k,k)$, but then $p\neq{\rm l.c.m.}(k,k)$.

If $X=S(2,3,3)$ then $d\myprod\ptil=12$.
Of course we cannot have $\ptil\dotsto 2$ and $\ptil\dotsto 3$. If $(\ptil,\ptil)\dotsto 2$
then $\ptil=2$ and $d=6$, so we get item 2 in~(\ref{Xtil:Spp:cand:list:eq}),
while if $(\ptil,\ptil)\dotsto 3$
then $\ptil=3$ and $d=4$, which is impossible since there would be
a partition of $4$ consisting of $3$'s only.
For $\ptil\dotsto 3$ and $\ptil\dotsto 3$ again $\ptil=3$ and $d=4$, whence
item 1 in~(\ref{Xtil:Spp:cand:list:eq}).

The discussion for $X=S(2,3,p)$ with $p=4,5$ is similar.
We get items 3 to 6 in~(\ref{Xtil:Spp:cand:list:eq}) for $p=4$ and
items 7 to 9 for $p=5$.
\end{proof}

\begin{prop}\label{Xtil:S22p:cand:list:prop}
The candidate surface branched covers
inducing some $S(2,2,\ptil)\dotsto X$ with $\ptil>1$ are
\begin{equation}\label{Xtil:S22p:cand:list:eq}
\begin{array}{ll}
S\argdotstoqui{4:1}{(2,1,1),(3,1),(4)}S &
S\argdotstoqui{6:1}{(2,2,1,1),(3,3),(4,2)}S \\
S\argdotstoqui{6:1}{(2,2,1,1),(3,3),(5,1)}S &
S\argdotstoqui{10:1}{(2,\ldots,2,1,1),(3,3,3,1),(5,5)}S \\
S\argdotstoqui{15:1}{(2,\ldots,2,1,1,1),(3,\ldots,3),(5,5,5)}S. &
\end{array}
\end{equation}
\end{prop}

\begin{proof}
Note first that $X$ cannot be $S$ or $S(p)$ since $d\geqslant 2$.
If $X=S(p,q)$ we have the following possibilities:
\begin{itemize}
\item $(2,2,\ptil)\dotsto p$. Then $p=2k$, so
$d=k+k+2k/\ptil+(m_1-3)  2k=m_2  q$;
\item $\ptil\dotsto p$ and $(2,2)\dotsto q$. Then $p=k \ptil$ and $q=2h$,
so $d=k+(m_1-1)  k \ptil=h+h+(m_2-2)  2h$;
\item $(2,\ptil)\dotsto p$ and $2\dotsto q$. Then
$p=2k$ and $q=2h$, so $d=k+2k/\ptil+(m_1-2)  2k=h+(m_2-1)  2h$.
\end{itemize}
Since $\frac1\ptil=\frac dp+\frac dq$, in all cases we deduce that $m_1+m_2=2$, which is absurd, so we do not get
any candidate cover.

Now suppose $X=S(2,2,p)$. We have the following possibilities:
$$\begin{array}{lll}
(2,2,\ptil)\dotsto 2 &
(2,2,\ptil)\dotsto p & \\
(2,2)\dotsto 2,\ \ptil\dotsto 2 &
(2,2)\dotsto 2,\ \ptil\dotsto p &
\ptil\dotsto 2,\ (2,2)\dotsto p \\
(2,\ptil)\dotsto 2,\ 2\dotsto 2 &
(2,\ptil)\dotsto 2,\ 2\dotsto p &
2\dotsto 2,\ (2,\ptil)\dotsto p,\\
2\dotsto 2,\ 2\dotsto 2,\ \ptil\dotsto p &
2\dotsto 2,\ \ptil\dotsto 2,\ 2\dotsto p
\end{array}$$
and $\ptil$ must actually be $2$ in items 1, 3, 5, 6, 7 and 10. Items 1, 3,
5, 6 and 8 are then impossible because $d$ should be both even and odd.
In item 2 we have $d=2k$ and $m_1=m_2=k$, whence $m_3=2$, which is impossible.
In items 4 and 7 we have $d=2k$, whence $m_1=k+1$ and $m_2=k$, so $m_3=1$, which is absurd.
In items 9 and 10 we have $d=2k+1$ and $m_1=m_2=k+1$, whence $m_3=1$, which is absurd.
This shows that there is no candidate cover for $X=S(2,2,p)$.

The cases where $X=S(2,3,p)$ for $p=3,4,5$ are easier to discuss and hence left to the reader.
For $p=3$ there is nothing, for $p=4$ there are items 1 and 2 in~(\ref{Xtil:S22p:cand:list:eq}),
and for $p=5$ items 3 to 5.
\end{proof}

The candidate surface branched covers inducing some
$X(2,3,\ptil)\dotsto X$ for $\ptil=3,4,5$ are not hard to analyze using
the same methods employed above,
so we will not spell out the proofs of the next two results.
The most delicate point is always to
exclude the cases $X=S(p,q)$ and $X=S(2,2,p)$.

\begin{prop}\label{Xtil:S233:cand:list:prop}
The only candidate surface branched cover inducing some
$S(2,3,3)\dotsto X$ is
$S\argdotstoter{5:1}{(2,2,1),(3,1,1),(5)}S$.
\end{prop}

\begin{prop}\label{Xtil:234/5:cand:list:prop}
There are no candidate surface branched covers inducing any
$S(2,3,4)\dotsto X$ or $S(2,3,5)\dotsto X$.
\end{prop}

\paragraph{Spherical structures}
As already mentioned, any good 2-orbifold $X$ with $\chiorb(X)>0$
has a spherical structure, given by the action of some
finite group $\Gamma$ of isometries on the metric sphere $\matS$.
We will now explicitly describe each relevant $\Gamma$, thus identifying
$X$ with the quotient $\matS/\Gamma$.
To this end we will always
regard $\matS$ as the unit sphere of $\matC\times\matR$.

\bigskip\noindent
\textit{The football}\quad The
geometry of $S(p,p)$ is very easy, even if
for consistency with what follows we will not give the easiest description.
Consider in $\matS$ a wedge with vertices at the poles $(0,\pm1)$ and
edges passing through $(1,0)$ and $({\rm e}^{i\pi/p},0)$, so the
width is $\pi/p$. Now define $\Gammatil_{(p,p)}$ as
the group of isometries of $\matS$ generated by the reflections in the
edges of the wedge, and $\Gamma_{(p,p)}$ as the its subgroup
of orientation-preserving isometries. Then $\Gamma_{(p,p)}$ is
generated by the rotation of angle $2\pi/p$ around the poles $(0,\pm 1)$,
a fundamental domain for $\Gamma_{(p,p)}$ is the union of any two wedges sharing an edge,
and $S(p,p)=\matS/\Gamma_{(p,p)}$. See Fig.~\ref{pp:22p:tess:fig}-left.
\begin{figure}
\centering
\input{gh_ellip1.pstex_t}
\mycap{Tessellations of $\matS$ by fundamental domains of $\Gammatil_{(p,p)}$ and
$\Gammatil_{(2,2,p)}$\label{pp:22p:tess:fig}}
\end{figure}

\bigskip\noindent
\textit{Triangular orbifolds}\quad
The remaining spherical 2-orbifolds $S(2,q,p)$
with either $q=2$ or $q=3$ and $p=3,4,5$ are called \emph{triangular}.
The corresponding group $\Gamma_{(2,q,p)}$ is the subgroup of orientation-preserving elements
of a group $\Gammatil_{(2,q,p)}$ generated by the reflections in the edges of
a triangle $\Delta_{(2,q,p)}$ with angles $\pi/2,\pi/q,\pi/p$.
A fundamental domain of $\Gamma_{(2,q,p)}$ is then the union of $\Delta_{(2,q,p)}$
with its image under any of the reflections in its edges.
If $A,B,C$ are the vertices of $\Delta_{(2,q,p)}$ then $\Gamma_{(2,q,p)}$
is generated by the rotations $\alpha,\beta,\gamma$ of angles $2\pi/2,2\pi/q,2\pi/p$
around $\pm A,\pm B,\pm C$, respectively, and a presentation of $\Gamma_{(2,q,p)}$
is given by
$$\Gamma_{(2,q,p)}=\langle\alpha,\beta,\gamma|\ \alpha^2=\beta^q=\gamma^p=\alpha\myprod\beta\myprod\gamma=1\rangle.$$
The main point here is of course that the triangles $\Delta_{(2,q,p)}$ exist in
$\matS$. The choice for $q=2$ and arbitrary $p$ is easy: $A=(1,0)$, $B=({\rm e}^{i\pi/p},0)$
and $C=(0,1)$, see Fig.~\ref{pp:22p:tess:fig}-right. For $q=3$ see
Figures~\ref{233:234:tess:fig}
\begin{figure}
\centering
\input{gh_ellip2.pstex_t}
\mycap{Tessellations of $\matS$ by triangular fundamental domains of
$\Gammatil_{(2,3,3)}$ and $\Gammatil_{(2,3,4)}$\label{233:234:tess:fig}}
\end{figure}
and~\ref{235:tess:fig}.
\begin{figure}
\centering
\input{gh_ellip3.pstex_t}
\mycap{Partial tessellation of $\matS$ by triangular fundamental domains of
$\Gammatil_{(2,3,5)}$\label{235:tess:fig}}
\end{figure}
In each case the pictures also show the images of $\Delta(2,q,p)$ under the action of
$\Gammatil_{(2,q,p)}$. A
fundamental domain for the group $\Gamma_{(2,q,p)}$
giving $S(2,q,p)$ is always the union of any two
triangles sharing an edge.
In all the pictures for $q=3$ we have $A=(1,0)$. Moreover for $p=3$
we have $B=\left(\sqrt{\frac 13},\sqrt{\frac 23}\right)$ and
$C=\left(\sqrt{\frac 13}+i\sqrt{\frac 23},0\right)$, while for $p=4$ we have
$B=\left(\sqrt{\frac 23},\sqrt{\frac 13}\right)$
and $C=\left(\sqrt{\frac 12}+i\sqrt{\frac 12},0\right)$; for $p=5$ the exact
values of the coordinates of $B$ and $C$ are more
complicated.

\begin{rem}
\emph{Each spherical orbifold $\Gammatil_{(2,3,p)}$ is
the symmetry group of a regular polyhedron
(a Platonic solid)
inscribed in $\matS$:
the tetrahedron for $p=3$, the octahedron (or its dual cube) for $p=4$,
and the dodecahedron (or its dual icosahedron) for $p=5$.}
\end{rem}

\paragraph{Geometric realizations}
Let a candidate orbifold cover $\Xtil\dotsto X$ with
covering instructions be induced by some candidate surface branched cover.
When $\Xtil$ is bad and $X$ is spherical, Lemma~\ref{no:bad:on:good:lem}
implies that $\Xtil\dotsto X$ is not realizable.
The results established above show that $X$ is never bad, so we are left
to show that the cover is realizable when both $\Xtil$ and $X$ are spherical.
The case $\Xtil=S$ is however easy, since a realization
is given precisely by the geometric universal cover
$\matS\to X$ we have fixed. We will now deal with all the other cases,
by constructing an isometry $\ftil:\matS\to\matS$ that induces a realization
$f:\Xtil\to X$ such that $\pi\compo\ftil=f\compo\pitil$.

\begin{prop}\label{Xtil:Spp:all:real:prop}
All the candidate surface branched covers of
Proposition~\ref{Xtil:Spp:cand:list:prop} are realizable.
\end{prop}

\begin{proof}
The following list describes the
candidate orbifold covers $\Xtil\dotsto X$ and the covering instructions associated to
the items in~(\ref{Xtil:Spp:cand:list:eq}), together with the corresponding isometry
$\ftil:\matS\to\matS$ inducing the desired $f$
via the geometric universal covers
$\pitil:\matS\to\Xtil$ and $\pi:\matS\to X$ fixed above.
\begin{itemize}

\item $S(3,3)\argdotstoter{4:1}{3\dotsto 3,\ 3\dotsto 3}S(2,3,3)$ and
$\ftil$ is the rotation around $(\pm i,0)$ sending $(0,1)$ to point
$B$ of Fig.~\ref{233:234:tess:fig}-left;

\item $S(2,2)\argdotstobis{6:1}{(2,2)\dotsto 2}S(2,3,3)$
and $\ftil={\rm rot}_{(\pm i,0)}(\pi/2)$;

\item $S(4,4)\argdotstobis{6:1}{(4,4) \dotsto 4}S(2,3,4)$ and $\ftil$ is
the identity;

\item $S(3,3)\argdotstobis{8:1}{(3,3) \dotsto 3}S(2,3,4)$ and $\ftil$ is
the rotation around $(\pm i, 0)$ sending $(0,1)$ to point
$B$ of Fig.~\ref{233:234:tess:fig}-right;

\item $S(2,2)\argdotstobis{12:1}{(2,2) \dotsto 2}S(2,3,4)$ and
$\ftil={\rm rot}_{(\pm i,0)}(\pi/2)$;

\item $S(2,2)\argdotstobis{12:1}{(2,2) \dotsto 4}S(2,3,4)$ and $\ftil$ is
the identity;

\item $S(5,5)\argdotstobis{12:1}{(5,5) \dotsto 5}S(2,3,5)$ and $\ftil$ is the
rotation around $(\pm i, 0)$ mapping $(0,1)$ to the point
labelled $*$ in Fig.~\ref{55:to:235:tess:fig};
\begin{figure}
\centering
\input{gh_55on235.pstex_t}
\mycap{Left and right, black: the tessellations for $\Gammatil_{(5,5)}$
and $\Gammatil_{(2,3,5)}$ as in Figures~\ref{pp:22p:tess:fig}-left and~\ref{235:tess:fig}.
Right, red: image of the tessellation on the left
under the isometry $\ftil:\matS\to\matS$ inducing a realization of
$S(5,5)\dotsto S(2,3,5)$. The picture on the right is partial\label{55:to:235:tess:fig}}
\end{figure}

\item $S(3,3)\argdotstobis{20:1}{(3,3) \dotsto 3}S(2,3,5)$ and $\ftil$ is
the rotation around $(\pm i, 0)$ mapping $(0,1)$ to the point $B$
in Fig.~\ref{235:tess:fig}.
See also Fig.~\ref{33:to:235:tess:fig};
\begin{figure}
\centering
\input{gh_33on235.pstex_t}
\mycap{Left and right, black: the tessellations for $\Gammatil_{(3,3)}$
and $\Gammatil_{(2,3,5)}$ as in Figures~\ref{pp:22p:tess:fig}-left
and~\ref{235:tess:fig}. Right, red: image of the tessellation on the left
under the isometry $\ftil:\matS\to\matS$ inducing a realization of
$S(3,3)\dotsto S(2,3,5)$. The picture on the right is partial\label{33:to:235:tess:fig}}
\end{figure}

\item $S(2,2)\argdotstobis{30:1}{(2,2) \dotsto 2}S(2,3,5)$ and $\ftil$
is the identity;

\item $S(2,2)\argdotstoter{2k+1:1}{2 \dotsto 2,\ 2 \dotsto 2}S(2,2,2k+1)$
and $\ftil={\rm rot}_{(\pm i,0)}(\pi/2)$;

\item $S(2,2)\argdotstobis{2k+2:1}{(2,2) \dotsto 2}S(2,2,2k+2)$ and
$\ftil={\rm rot}_{(\pm i,0)}(\pi/2)$.
\end{itemize}
The proof is complete.
\end{proof}

\begin{prop}\label{Xtil:S22p:all:real:prop}
All the candidate surface branched covers of
Proposition~\ref{Xtil:S22p:cand:list:prop} are realizable.
\end{prop}

\begin{proof}
The candidate orbifold covers with covering instructions associated to
the items in~(\ref{Xtil:S22p:cand:list:eq}), and the corresponding isometries
$\ftil:\matS\to\matS$, are as follows:

\begin{itemize}

\item $S(2,2,3)\argdotstoter{4:1}{(2,2) \dotsto 2,\ 3 \dotsto 3}S(2,3,4)$ and $\ftil$ is
the rotation around $(\pm 1,0)$ mapping $(0,1)$ to point
$B$ in Fig.~\ref{233:234:tess:fig}-right followed by a rotation of angle $\pi/6$
around $\pm B$.
Fig.~\ref{223:to:234:tess:planar:fig};
\begin{figure}
\centering
\input{gh_223on234.pstex_t}
\mycap{Black: planar stereographic projection of
the tessellation for $\Gammatil_{(2,3,4)}$ from Fig.~\ref{233:234:tess:fig}-right.
Red: image of the tessellation for
$\Gammatil_{(2,2,3)}$ from Fig.~\ref{pp:22p:tess:fig}-right
under the isometry $\ftil:\matS\to\matS$ inducing a realization of
$S(2,2,3)\dotsto S(2,3,4)$\label{223:to:234:tess:planar:fig}}
\end{figure}

\item $S(2,2,2)\argdotstoter{6:1}{(2,2) \dotsto 2,\ 2 \dotsto 4}S(2,3,4)$ and $\ftil$ is
the identity;

\item $S(2,2,5)\argdotstoter{6:1}{(2,2) \dotsto 2,\ 5 \dotsto 5}S(2,3,5)$ and $\ftil$ is
the rotation around $(\pm 1, 0)$,
sending $(0,1)$ to the point of order $5$ best
visible in the upper hemisphere in Fig.~\ref{235:tess:fig}.
See also Fig.~\ref{225:to:235:tess:fig};
\begin{figure}
\centering
\input{gh_225on235.pstex_t}
\mycap{A planar illustration of the isometry $\ftil:\matS\to\matS$
inducing a realization of the cover $S(2,2,5)\dotsto S(2,3,5)$, with colours as above\label{225:to:235:tess:fig}}
\end{figure}

\item $S(2,2,3)\argdotstoter{10:1}{(2,2) \dotsto 2,\ 3 \dotsto 3}S(2,3,5)$ and $\ftil$ is
the rotation around $(\pm i, 0)$ that maps the pole $(0,1)$ to point $B$
in Fig.~\ref{235:tess:fig} followed by a rotation of angle $\pi/6$ around $\pm B$.
See also
Fig.~\ref{223:to:235:pl:tess:fig};
\begin{figure}
\centering
\input{gh_223on235.pstex_t}
\mycap{A planar illustration of the isometry $\ftil$ inducing a realization
of the cover $S(2,2,3)\dotsto S(2,3,5)$, with colours as above\label{223:to:235:pl:tess:fig}}
\end{figure}

\item $S(2,2,2)\argdotstoter{15:1}{(2,2,2) \dotsto 2}S(2,3,5)$ and $\ftil$ is
the identity.
\end{itemize}
The proof is complete.
\end{proof}

\begin{prop}\label{Xtil:S233:real:prop}
The candidate surface branched cover of
Proposition~\ref{Xtil:S233:cand:list:prop}
is realizable.
\end{prop}

\begin{proof}
The candidate orbifold cover is in this case
$$S(2,3,3)\argdotstoter{5:1}{2\dotsto 2,\ (3,3)\dotsto 3}S(2,3,5)$$
and the isometry $\ftil:\matS\to\matS$ inducing its realization
is the rotation of angle $\pi/4$ around $(\pm1,0)$. See also
Fig.~\ref{233:to:235:tess:fig}.
\begin{figure}
\centering
\input{gh_233on235.pstex_t}
\mycap{A partial planar illustration of the isometry $\ftil:\matS\to\matS$ inducing a realization
of the cover $S(2,3,3)\dotsto S(2,3,5)$, with colours as
above\label{233:to:235:tess:fig}}
\end{figure}
\end{proof}

\section{The Euclidean case}\label{zero:chi:sec}
In this section we investigate realizability of candidate surface branched covers inducing
candidate orbifold covers $\Xtil\dotsto X$ with $\chiorb(\Xtil)=\chiorb(X)=0$.
This means that $\Xtil$ and $X$ must belong to the list
$$T,\quad S(2,4,4),\quad S(2,3,6),\quad S(3,3,3),\quad S(2,2,2,2)$$
where $T$ is the torus. Recalling that the orders
of the cone points of $\Xtil$ must divide those of $X$, we see that the only
possibilities are the cases (0) to (7) shown in Fig.~\ref{gen:eucl:fig},
\begin{figure}
\centering
\input{gh_geneucl.pstex_t}
\mycap{Possible covers between Euclidean orbifolds\label{gen:eucl:fig}}
\end{figure}
that we will analyze using Euclidean geometry. Namely:
\begin{itemize}
\item We will fix on $X$ a Euclidean structure given by some
$\pi:\matE\to X$;
\item We will assume that
$\Xtil\dotsto^{d:1} X$ is realized by some
map $f$, we will use Lemma~\ref{Eucl:X:then:Eucl:Xtil:prop}
to deduce there is a corresponding affine map $\ftil:\matE\to\matE$,
and we will analyze $\ftil$ to show that $d$ must satisfy certain conditions;
\item We will employ the calculations of the previous
point to show that if $d$ satisfies the conditions
then $\ftil$, whence $f$, exists.
\end{itemize}

\paragraph{Euclidean structures}
On each of the three triangular orbifolds $S(p,q,r)$
with $\frac1p+\frac1q+\frac1r=1$ the Euclidean structure
(unique up to rescaling) is constructed essentially as in the
spherical case. We take a triangle $\Delta(p,q,r)$ in $\matE$
with angles $\pi/p,\pi/q,\pi/r$, the group $\Gammatil_{(p,q,r)}$ generated by
the reflections in the edges of $\Delta(p,q,r)$, and its subgroup
$\Gamma_{(p,q,r)}$ of orientation-preserving isometries.
Then $S(p,q,r)$ is the quotient of $\matE$ under the action of
$\Gamma_{(p,q,r)}$, which is generated by the rotations
of angles $2\pi/p,2\pi/q,2\pi/r$ around the vertices of
$\Delta(p,q,r)$, and a fundamental
domain is given by the union of $\Delta(p,q,r)$ with any of its reflected copies
in one of the edges. In each case we now make a precise choice
for the vertices $\Atil^{(p)},\Btil^{(q)},\Ctil^{(r)}$
of $\Delta(p,q,r)$ and determine
the resulting area $\calA$ of $S(p,q,r)$. See also
Fig.~\ref{fundomains:fig}.
\begin{figure}
\centering
\input{gh_fundom1.pstex_t}
\mycap{The fixed fundamental domains for $S(2,4,4)$, $S(2,3,6)$ and $S(3,3,3)$\label{fundomains:fig}}
\end{figure}
$$\begin{array}{lllll}
\Delta(2,4,4): & \Atil^{(2)}=1, &
\Btil^{(4)}=0, & \Ctil^{(4)}=1+i, & \calA(S(2,4,4))=1,\\
\Delta(2,3,6): & \Atil^{(2)}=\frac12, &
\Btil^{(3)}=0, & \Ctil^{(6)}=\frac{1+i\sqrt3}2, & \calA(S(2,3,6))=\frac{\sqrt{3}}4,\\
\Delta(3,3,3): & \Atil^{(3)}=0, &
\Btil^{(3)}=1, & \Ctil^{(3)}=\frac{1+i\sqrt3}2, & \calA(S(3,3,3))=\frac{\sqrt{3}}2.
\end{array}$$

The situation for $S(2,2,2,2)$ is slightly different, since there is flexibility besides rescaling.
For $s,t\in\matR$ with $s>0$ we consider in $\matE$ the quadrilateral $Q_{s,t}$ with corners
$$\Atil^{(2)}=0,\qquad \Btil^{(2)}=\frac1s+it,\qquad \Ctil^{(2)}=\frac1s+i(s+t),\qquad \Dtil^{(2)}=is$$
and we define $\Gamma^{s,t}_{(2,2,2,2)}$ as the group generated by the rotations of angle $\pi$ around these
points. Then the action of
$\Gamma^{s,t}_{(2,2,2,2)}$ on $\matE$ defines on $S(2,2,2,2)$ a Euclidean
structure of area $2$, and a fundamental domain is
given by the union of $Q_{s,t}$ with any translate of itself having an edge in common
with itself, as shown in
Fig.~\ref{fundomains:bis:fig}-left.
\begin{figure}
\centering
\input{gh_fundom2.pstex_t}
\mycap{The fundamental domains for $S(2,2,2,2)$ for general $s,t$ and for $s=t=1$\label{fundomains:bis:fig}}
\end{figure}
When $S(2,2,2,2)$ plays the r\^ole of $X$ in $\Xtil\dotsto X$ we will
endow it with the structure given by $s=t=1$, as shown in
Fig.~\ref{fundomains:bis:fig}-right. It is an easy exercise to
check that any other structure with area $2$ is defined by $\Gamma^{s,t}_{(2,2,2,2)}$ for
some $s,t$.

\paragraph{General geometric tools}
Our analysis of the candidate covers (1)-(7) of Fig.~\ref{gen:eucl:fig}
relies on certain facts that we will use repeatedly. The first is the exact determination of
the lifts of the cone points, that we now describe.
For any of our four Euclidean $X$'s, with the structure $\pi:\matE\to X=\matE/\Gamma$ we have fixed,
and any vertex $\widetilde V^{(p)}$ of the fundamental domain for $\Gamma$ described above, we set
$V^{(p)}=\pi(\widetilde V^{(p)})$, so that its cone order is $p$. Then $\pi^{-1}(V^{(p)})$
will be some set $\{\widetilde V^{(p)}_j\}\subset\matE$, with $j$ varying in a suitable set of indices.
The exact lists of lifts are as follows:

\begin{equation}\label{244:lifts:eq}
\begin{array}{ll}
S(2,4,4): & ({\hbox{\rm see\ Fig.~\ref{244:2222:tess:fig}-left}})\\
\Atil^{(2)}_{a,b}=a+ib\quad & a,b\in\matZ,\ a\nequiv b\ ({\rm mod}\ 2)\\
\Btil^{(4)}_{a,b}=a+ib\quad & a,b\in\matZ,\ a\equiv b\equiv 0\ ({\rm mod}\ 2)\\
\Ctil^{(4)}_{a,b}=a+ib\quad & a,b\in\matZ,\ a\equiv b\equiv 1\ ({\rm mod}\ 2);
\end{array}
\end{equation}
\begin{figure}
\centering
\input{gh_eucl_tess1.pstex_t}
\mycap{Tessellations of $\matE$ induced by the geometric structures
fixed on $S(2,4,4)$ and $S(2,2,2,2)$.\label{244:2222:tess:fig}}
\end{figure}

\begin{equation}\label{2222:lifts:eq}
\begin{array}{ll}
S(2,2,2,2): & ({\hbox{\rm see\ Fig.~\ref{244:2222:tess:fig}-right}})\\
\Atil^{(2)}_{a,b}=a+ib\quad & a,b\in\matZ,\ a\equiv b\equiv 0\ ({\rm mod}\ 2) \\
\Btil^{(2)}_{a,b}=a+ib\quad & a,b\in\matZ,\ a\equiv 1,\ b\equiv 0\ ({\rm mod}\ 2) \\
\Ctil^{(2)}_{a,b}=a+ib\quad & a,b\in\matZ,\ a\equiv b\equiv 1\ ({\rm mod}\ 2) \\
\Dtil^{(2)}_{a,b}=a+ib\quad & a,b\in\matZ,\ a\equiv 0,\ b\equiv 1\ ({\rm mod}\ 2);
\end{array}
\end{equation}

\begin{equation}\label{236:lifts:eq}
\begin{array}{ll}
S(2,3,6): & {\rm with}\ \omega=\frac{1+i\sqrt3}2\ ({\hbox{\rm see\ Fig.~\ref{236:333:tess:fig}-left}})\\
\Atil^{(2)}_{a,b}=\frac{1}{2}(a+\omega b)\quad & a,b\in\matZ\ {\rm not\ both\ even},\ a-b\equiv 1\ ({\rm mod}\ 3)\\
\Btil^{(3)}_{a,b}=a+\omega b\quad & a,b\in\matZ,\ a-b\nequiv 2\ ({\rm mod}\ 3)\\
\Ctil^{(6)}_{a,b}=a+\omega b\quad & a,b\in\matZ,\ a-b\equiv 2\ ({\rm mod}\ 3);
\end{array}
\end{equation}
\begin{figure}
\centering
\input{gh_eucl_tess2.pstex_t}
\mycap{Tessellations of $\matE$ induced by the geometric structure
fixed on $S(2,3,6)$ and $S(3,3,3)$\label{236:333:tess:fig}}
\end{figure}

\begin{equation}\label{333:lifts:eq}
\begin{array}{ll}
S(3,3,3): & {\rm with}\ \omega=\frac{1+i\sqrt3}2\ ({\hbox{\rm see\ Fig.~\ref{236:333:tess:fig}-right}})\\
\Atil^{(3)}_{a,b}=a+\omega b,\quad & a,b\in\matZ,\ a-b\equiv 0\ ({\rm mod}\ 3)\\
\Btil^{(3)}_{a,b}=a+\omega b,\quad & a,b\in\matZ,\ a-b\equiv 1\ ({\rm mod}\ 3)\\
\Ctil^{(3)}_{a,b}=a+\omega b,\quad & a,b\in\matZ,\ a-b\equiv 2\ ({\rm mod}\ 3).
\end{array}
\end{equation}

\bigskip

Also important for us will be the symmetries of our Euclidean $X$'s. More precisely
we will use the following facts:
\begin{itemize}
\item There exists a symmetry of $S(2,4,4)$ switching $B^{(4)}$ and $C^{(4)}$;
\item Any permutation of $A^{(3)},B^{(3)},C^{(3)}$ is induced by a symmetry of $S(3,3,3)$;
\item Any permutation of $A^{(2)},B^{(2)},C^{(2)},D^{(2)}$ is induced by a
symmetry of $S(2,2,2,2)$.
\end{itemize}
It is perhaps worth noting that all the symmetries for $S(2,4,4)$ and $S(3,3,3)$ are isometries
with respect to the structures we have fixed, and those of $S(2,2,2,2)$ are isometries
if we endow it with the structure such that the triangle with vertices
$\Atil^{(2)},\Btil^{(2)},\Dtil^{(2)}$ is equilateral.

\bigskip

Here is another tool we will use very often.
As one sees, if $\Gamma<\Isomp(\matE)$
defines a Euclidean structure on a 2-orbifold $X=\matE/\Gamma$
then $\Gamma$ has a maximal torsion-free subgroup, a rank-2 lattice $\Lambda(\Gamma)$.
Denoting for $u\in\matC$ by $\tau_u$
the translation $z\mapsto z+u$, the lattices for the groups we have fixed
are as follows:
$$\begin{array}{ll}
\Lambda_{(2,4,4)}=\langle \tau_2,\tau_{2i}\rangle,\qquad
\Lambda^{s,t}_{(2,2,2,2)}=\left\langle\tau_{2is},\tau_{2\left(\frac{1}{s}+it\right)}\right\rangle,\\
\Lambda_{(2,3,6)}=\Lambda_{(3,3,3)}=\left\langle \tau_{i\sqrt3},\tau_{\frac{3+i\sqrt3}2}\right\rangle.
\end{array}$$
Moreover the following holds:
\begin{lem}\label{lattice:to:lattice:lem}
Let $\Gammatil,\Gamma<\Isomp(\matE)$ define Euclidean orbifolds
$\Xtil=\matE/\Gammatil$ and $X=\matE/\Gamma$.
Suppose that $\Lambda(\Gammatil)=\langle \tau_{\util_1},\tau_{\util_2}\rangle$
and $\Lambda(\Gamma)=\langle \tau_{u_1},\tau_{u_2}\rangle$.
Let $\ftil:\matE\to\matE$ given by $\ftil(z)=\lambda\myprod z+\mu$
induce an orbifold cover $\Xtil\to X$. Then
$\lambda\myprod\util_1$ and $\lambda\myprod\util_2$ are integer
linear combinations of $u_1$ and $u_2$.
\end{lem}

\begin{proof}
The map $\ftil$ induces a homomorphism $\ftil_*:\Gammatil\to\Gamma$
given by $\ftil_*(\tau_u)=\tau_{\lambda\cdot u}$ which maps
$\Lambda(\Gammatil)$ to $\Lambda(\Gamma)$,
and the conclusion easily follows.
\end{proof}

\paragraph{Exceptions and geometric realizations}
We will now state and prove 8 theorems corresponding to
the cases (0)-(7) of Fig.~\ref{gen:eucl:fig},
thus establishing Theorem~\ref{zero:chi:main:thm}.
Cases (1)-(3) imply in particular Theorems~\ref{244:main:thm},~\ref{236:main:thm},
and~\ref{333:main:thm}.

\begin{thm}[case (0) in Fig.~\ref{gen:eucl:fig}]\label{T:on:X:thm}
The candidate surface branched covers inducing some
$T\dotsto X$ are
$$\begin{array}{lll}
T\dotsto^{k:1}T &
T\argdotstosex{2k:1}{(2,\ldots,2),(2,\ldots,2),(2,\ldots,2),(2,\ldots,2)}S\ &
T\argdotstoqui{3k:1}{(3,\ldots,3),(3,\ldots,3),(3,\ldots,3)}S\\
T\argdotstoqui{4k:1}{(2,\ldots,2),(4,\ldots,4),(4,\ldots,4)}S\ &
T\argdotstosex{6k:1}{(2,\ldots,2),(3,\ldots,3),(6,\ldots,6)}S &
\end{array}$$
with arbitrary $k\geqslant 1$, and they are all realizable.
\end{thm}

\begin{proof}
The first assertion and realizability of any $T\dotsto^{k:1}T$ are easy.
For any $X\neq T$ let $X=\matE/\Gamma$ and identify $T$ to
$\matE/\Lambda(\Gamma)$. Since $\Lambda(\Gamma)<\Gamma$ we have
an induced orbifold cover $T\to X$, which realizes the relevant
$T\dotsto X$ in the special case $k=1$. The conclusion follows by taking compositions.
\end{proof}

\begin{thm}[case (1) in Fig.~\ref{gen:eucl:fig}]\label{244:to:244:thm}
The candidate surface branched covers inducing some
$S(2,4,4)\dotsto^{d:1} S(2,4,4)$ are
$$
S\argdotstoqui{4k+1:1}{(2,\dots,2,1),(4,\dots,4,1),(4,\dots,4,1)}S\quad
S\argdotstoqui{4k+2:1}{(2,\dots,2),(4,\dots,4,2),(4,\dots,4,1,1)}S\quad
S\argdotstoqui{4k+4:1}{(2,\dots,2),(4,\dots,4),(4,\dots,4,2,1,1)}S
$$
for $k\geqslant 1$, and they are realizable if and only if, respectively:
\begin{itemize}
\item $d=x^2+y^2$ for some $x,y\in\matN$ of different parity;
\item $d=2(x^2+y^2)$ for some $x,y\in\matN$ of different parity;
\item $d=4(x^2+y^2)$ for some $x,y\in\matN$.
\end {itemize}
\end{thm}

\begin{proof}
A candidate $S(2,4,4)\dotsto^{d:1} S(2,4,4)$
can be complemented with the covering instructions
$$\begin{array}{lll}
2\dotsto2,\ 4\dotsto4,\ 4\dotsto4,\quad &
2\dotsto4,\ (4,4)\dotsto4,\quad &
(2,4,4)\dotsto 4,\quad\\
2\dotsto2,\ (4,4)\dotsto 4,\quad &
(2,4)\dotsto 4,\ 4\dotsto 4\quad &
\end{array}$$
and it is very easy to see that the first three come from the
candidate surface branched covers of the statement, while the last two
do not come from any candidate cover (recall that the simplified version~(\ref{RH:SS3:eq})
of the Riemann-Hurwitz formula must be satisfied).

Now suppose there is a realization $f:S(2,4,4)\myto^{d:1} S(2,4,4)$ of one of the three relevant
candidate covers. Proposition~\ref{Eucl:X:then:Eucl:Xtil:prop} implies that there is
a geometric universal cover $\pitil:\matE\to S(2,4,4)$
and an affine map $\ftil:\matE\to\matE$ with $\pi\compo\ftil=\pitil\compo f$. But
the Euclidean structure of $S(2,4,4)$ is unique up to scaling, so $\pitil=\pi$.
If $\ftil(z)=\lambda\myprod z+\mu$, since
$\Lambda_{(2,4,4)}=\langle \tau_2,\tau_{2i}\rangle$,
Lemma~\ref{lattice:to:lattice:lem} implies that $\lambda=n+im$ for some
$n,m\in\matZ$, whence $d=n^2+m^2$.

\medskip

We now employ the notation of~(\ref{244:lifts:eq}) and
note that for all three candidates we can assume, by the symmetry of $S(2,4,4)$,
that $f\left(B^{(4)}\right)=B^{(4)}$, whence that
$\ftil\left(\Btil^{(4)}_{0,0}\right)=\Btil^{(4)}_{0,0}$, namely $\mu=0$.
We then proceed separately for the three candidates. For the first one
we have $f\left(A^{(2)}\right)=A^{(2)}$
so $\ftil\left(\Atil^{(2)}_{1,0}\right)=\lambda$ is some $\Atil^{(2)}_*$. Therefore
$n$ and $m$ have different parity, and we can set $x=|n|$ and $y=|m|$ getting that
$d=x^2+y^2$ for $x,y\in\matN$ of different parity. Conversely if $d$ has this form we
define $\ftil(z)=(x+iy)\myprod z$. Then
$$\ftil\left(\Atil^{(2)}_{1,0}\right)=\Atil^{(2)}_{x,y},\quad
\ftil\left(\Btil^{(4)}_{0,0}\right)=\Btil^{(4)}_{0,0},\quad
\ftil\left(\Ctil^{(4)}_{1,1}\right)=(x-y)+i(x+y)=\Ctil^{(4)}_{x-y,x+y}$$
where the last equality depends on the fact that $x-y\equiv x+y\equiv 1\ ({\rm mod}\ 2)$.
It easily follows that $\ftil$ induces a realization of the candidate.

\medskip

For the second candidate
$f\left(A^{(2)}\right)=C^{(4)}$,
hence $\ftil\left(\Atil^{(2)}_{1,0}\right)=\lambda$ is some
$\Ctil^{(4)}_*$, namely $n$ and $m$ are odd. Setting
$x=\frac12|n+m|$ and $y=\frac12|n-m|$ we
see that $x,y\in\matN$ have different parity and $d=2(x^2+y^2)$.
Conversely if $d=2(x^2+y^2)$ with $x,y$ of different parity, we define
$\ftil(z)=((x+y)+i(x-y))\myprod z$. Then
$$\ftil\left(\Atil^{(2)}_{1,0}\right)=\Ctil^{(4)}_{x+y,x-y},\quad
\ftil\left(\Btil^{(4)}_{0,0}\right)=\Btil^{(4)}_{0,0},\quad
\ftil\left(\Ctil^{(4)}_{1,1}\right)=\Btil^{(4)}_{2y,2x}$$
from which it is easy to see that $\ftil$
induces a realization of the candidate.

\medskip

For the last candidate
$f\left(A^{(2)}\right)=B^{(4)}$
hence $\ftil\left(\Atil^{(2)}_{1,0}\right)$ is some $\Btil^{(4)}_*$,
so $n$ and $m$ are even. Setting $x=\frac12|n|$ and $y=\frac12|m|$
we see that $d=4(x^2+y^2)$ for $x,y\in\matN$.
Conversely if $d=4(x^2+y^2)$ and we define
$\ftil(z)=(2x+2iy)\myprod z$ then
$$\ftil\left(\Atil^{(2)}_{1,0}\right)=\Btil^{(4)}_{2x,2y},\quad
\ftil\left(\Btil^{(4)}_{0,0}\right)=\Btil^{(4)}_{0,0},\quad
\ftil\left(\Ctil^{(4)}_{1,1}\right)=\Btil^{(4)}_{2(x-y),2(x+y)}$$
therefore $\ftil$ induces a realization of the candidate.
\end{proof}

\begin{rem}\label{two:thirds:rem}
\emph{One feature of the above proof is worth pointing out.
After assuming that some $S(2,4,4)\dotsto^{d:1} S(2,4,4)$ is realized
by some map, we have always used only ``two thirds'' of the branching instruction to
show that $d$ has the appropriate form. The same phenomenon will occur in
all the next proofs, except that of Theorem~\ref{333:to:236:thm}.
Note however that to check that some $\ftil$ defined starting from a degree
$d$ with appropriate form induces a realization of the corresponding candidate,
we had to check all three conditions.}
\end{rem}

As an illustration of the previous proof, we provide in Fig.~\ref{244:on:244:fig}
\begin{figure}
\centering
\input{gh_244on244.pstex_t}
\mycap{A map $\ftil:\matE\to\matE$ inducing a degree-5 cover
of $S(2,4,4)$ onto itself\label{244:on:244:fig}}
\end{figure}
a description of the map $\ftil:\matE\to\matE$ inducing a realization
of $S\argdotstoter{5:1}{(2,2,1),(4,1),(4,1)}S$.

\begin{thm}[case (2) in Fig.~\ref{gen:eucl:fig}]\label{236:to:236:thm}
The candidate surface branched covers inducing some $S(2,3,6)\dotsto^{d:1} S(2,3,6)$ are
$$\begin{array}{ll}
S\argdotstosex{6k+1:1}{(2,\dots,2,1),(3,\dots,3,1),(6,\dots,6,1)}S\quad &
S\argdotstosex{6k+3:1}{(2,\dots,2,1),(3,\dots,3),(6,\dots,6,2,1)}S\quad\\
S\argdotstosex{6k+4:1}{(2,\dots,2),(3,\dots,3,1),(6,\dots,6,3,1)}S\quad &
S\argdotstosex{6k+6:1}{(2,\dots,2),(3,\dots,3),(6,\dots,6,3,2,1)}S\quad
\end{array}$$
with $k\geqslant 1$, and they are realizable if and only if, respectively:
\begin{itemize}
\item $d=x^2+xy+y^2$ with $x,y\in\matN$ not both even and $x\nequiv y\ ({\rm mod}\ 3)$;
\item $d=3(x^2+3xy+3y^2)$ with $x,y\in\matN$ not both even;
\item $d=12(x^2+3xy+3y^2)+16$ with $x,y\in\matN$;
\item $d=12(x^2+3xy+3y^2)$ with $x,y\in\matN$.
\end{itemize}
\end{thm}

\begin{proof}
A candidate orbifold cover $S(2,3,6)\dotsto S(2,3,6)$ can be complemented with
the covering instructions
$$\begin{array}{ll}
2\dotsto 2,\ 3\dotsto 3,\ 6\dotsto 6\quad &
2\dotsto 2,\ (3,6)\dotsto 6,\quad\\
3\dotsto 3,\ (2,6)\dotsto 6,\quad &
(2,3,6)\dotsto 6
\end{array}$$
which are induced by the candidate surface branched covers of the statement
(formula~(\ref{RH:SS3:eq}) is always satisfied in this case).

The scheme of the proof is now as for case (1): we consider the universal cover $\pi:\matE\to S(2,3,6)$
we have fixed, we assume that a map $f:S(2,3,6)\myto^{d:1} S(2,3,6)$
realizing some candidate cover exists, we use Proposition~\ref{Eucl:X:then:Eucl:Xtil:prop}
to find $\ftil:\matE\to\matE$ with $\ftil(z)=\lambda\myprod z+\mu$ and
$\pi\compo\ftil=f\compo\pi$, and we show that $d=|\lambda|^2$ has the
appropriate form. Moreover essentially the same calculations
will also allow us to prove the converse.
We will always use the notation fixed
in~(\ref{236:lifts:eq}), in particular $\omega=\frac{1+i\sqrt3}2$.

We first apply Lemma~\ref{lattice:to:lattice:lem}. Since $\Lambda_{(2,3,6)}=
\left\langle \tau_{i\sqrt3},\tau_{\frac{3+i\sqrt3}2}\right\rangle$ we must have
$$\lambda\myprod i\sqrt3=n\myprod i\sqrt3+m\myprod\frac{3+i\sqrt3}2$$
for some $n,m\in\matZ$, which easily implies that
$$\lambda=n+\frac m2-i\sqrt{3}\frac m2=(n+m)-\omega m$$
therefore $d=n^2+nm+m^2$. Notice that Lemma~\ref{lattice:to:lattice:lem}
does not give any other condition, because
$\lambda\myprod\frac{3+i\sqrt3}2=(-m)\myprod i\sqrt 3+(n+m)\myprod \frac{3+i\sqrt3}2$.

\medskip

We now analyze our four candidates, starting from the first one. Since
$$f\left(A^{(2)}\right)=A^{(2)},\quad
f\left(B^{(3)}\right)=B^{(3)},\quad
f\left(C^{(6)}\right)=C^{(6)}$$
we can assume $\ftil\left(\Btil^{(3)}_{0,0}\right)=\Btil^{(3)}_{0,0}$, so
$\mu=0$, and
$\ftil\left(\Atil^{(2)}_{1,0}\right)=\Atil^{(2)}_{a,b}$
for some $a,b\in\matZ$ not both even with $a-b\equiv 1\ ({\rm mod}\ 3)$.
Since $a=n+m$ and $b=-m$ we deduce that $n,m$ are not both even and $n-m\equiv 1\ ({\rm mod}\ 3)$.
If $n,m\geqslant 0$ or $n,m\leqslant 0$ we set $x=|n|$ and $y=|m|$, getting
$d=x^2+xy+y^2$ with $x,y\in\matN$ not both even and $x\nequiv y\ ({\rm mod}\ 3)$.
Otherwise we have $n>0>m$ up to permutation. If $n\geqslant -m$ we set $x=n+m$ and $y=-m$, otherwise
$x=n$ and $y=-n-m$, and again we have
$d=x^2+xy+y^2$ with $x,y\in\matN$ not both even and $x\nequiv y\ ({\rm mod}\ 3)$.

Conversely, assume $d=x^2+xy+y^2$ with $x,y\in\matN$ not both even and
$x\nequiv y\ ({\rm mod}\ 3)$. Up to changing sign to both $x$ and $y$ we can suppose that
$x-y\equiv 1\ ({\rm mod}\ 3)$ and define $\ftil(z)=((x+y)-\omega y)\myprod z$.
Of course $\ftil\left(\Btil^{(3)}_{0,0}\right)=\Btil^{(3)}_{0,0}$,
the above calculations show that $\ftil\left(\Atil^{(2)}_{1,0}\right)=\Atil^{(2)}_{x+y,-y}$,
and using the identity $\omega^2=\omega-1$ we have
$\ftil\left(\Ctil^{(6)}_{0,1}\right)=
(x+y-\omega y)\myprod\omega=y+\omega x=\Ctil^{(6)}_{y,x}$
where the last equality depends on the fact that
$y-x\equiv -(x-y)\equiv 2\ ({\rm mod}\ 3)$.
This easily implies that $\ftil$ induces a realization of the first candidate.

\medskip

Let us turn to the second candidate.
Since
$$f\left(A^{(2)}\right)=A^{(2)},\quad
f\left(B^{(3)}\right)=C^{(6)},\quad
f\left(C^{(6)}\right)=C^{(6)}$$
we can assume $\ftil\left(\Btil^{(3)}_{0,0}\right)=\Ctil^{(6)}_{0,1}$,
namely $\mu=\omega$. In addition
$\ftil\left(\Atil^{(2)}_{1,0}\right)=\Atil^{(2)}_{a,b}$
for some $a,b\in\matZ$ not both even with $a-b\equiv 1\ ({\rm mod}\ 3)$.
Since $a=n+m$ and $b=2-m$ we deduce that $n,m$ are not both even and $n\equiv m\ ({\rm mod}\ 3)$.
Setting $x=n$ and $y=(m-n)/3$ we then get $d=3(x^2+3xy+3y^2)$ for $x,y\in\matZ$ not both even.
Reducing to the case $x,y\in\matN$ is a routine matter that we leave to the
reader.

Conversely, assume $d=3(x^2+3xy+3y^2)$ for $x,y\in\matN$ not both even,
set $n=x$ and $m=x+3y$ and define
$$\ftil(z)=((n+m)-\omega m)\myprod z+\omega.$$
Then $\ftil\left(\Atil^{(2)}_{1,0}\right)$ is some $\Atil^{(2)}_*$ by the above calculations, while
$$\ftil\left(\Ctil^{(6)}_{0,1}\right)=\ftil(\omega)=(x+3y)+\omega(x+1)$$
which is some $\Ctil^{(6)}_*$ because $(x+3y)-(x+1)\equiv 2\ ({\rm mod}\ 3)$.
It follows that $\ftil$ induces a realization of the candidate.

\medskip

For the third candidate we have
$$f\left(A^{(2)}\right)=C^{(6)},\quad
f\left(B^{(3)}\right)=B^{(3)},\quad
f\left(C^{(6)}\right)=C^{(6)}$$
so we can take $\mu=0$ and
$\ftil\left(\Atil^{(2)}_{1,0}\right)=\Ctil^{(6)}_{a,b}$
for some $a,b\in\matZ$ with $a-b\equiv 2\ ({\rm mod}\ 3)$.
Since $a=(n+m)/2$ and $b=-m/2$ we have that $n$ and $m$ are even and
$n-m\equiv 1\ ({\rm mod}\ 3)$. So we can define
$x=\frac{n+4}2$ and $y=\frac{m-n-8}6$ and
we have $d=12(x^2+3xy+3y^2)+16$ with $x,y\in\matZ$. Again it is easy to reduce to the case
$x,y\in\matN$.

Conversely, suppose $d=12(x^2+3xy+3y^2)+16$ with $x,y\in\matN$, set
$n=2(x-2)$ and $m=2(x+3y+2)$, and define
$\ftil(z)=((n+m)-\omega m)\myprod z$. Then
of course $\ftil\left(\Btil^{(3)}_{0,0}\right)=\Btil^{(3)}_{0,0}$,
the above calculations show that $\ftil\left(\Atil^{(2)}_{1,0}\right)$ is some $\Ctil^{(6)}_*$, while
$$\ftil\left(\Ctil^{(6)}_{0,1}\right)=\ftil(\omega)=2(x+3y+2)+\omega 2(x-2)$$
which is some $\Ctil^{(6)}_*$ because $2(x+3y+2)-2(x-2)\equiv 2\ ({\rm mod}\ 3)$,
hence $\ftil$ induces a realization of the candidate.

\medskip

For the last candidate
$$f\left(A^{(2)}\right)=C^{(6)},\quad
f\left(B^{(3)}\right)=C^{(6)},\quad
f\left(C^{(6)}\right)=C^{(6)}$$
so we can assume $\ftil\left(\Btil^{(3)}_{0,0}\right)=\Ctil^{(6)}_{0,1}$,
namely $\mu=\omega$. In addition
$\ftil\left(\Atil^{(2)}_{1,0}\right)=\Ctil^{(6)}_{a,b}$
for some $a,b\in\matZ$ with $a-b\equiv 2\ ({\rm mod}\ 3)$.
Since $a=(n+m)/2$ and $b=(2-m)/2$ we deduce that $n,m$ are both even and $n\equiv m\ ({\rm mod}\ 3)$.
Setting $x=n/2$ and $y=(m-n)/6$ we then get $d=12(x^2+3xy+3y^2)$ for $x,y\in\matZ$, and again
we can reduce to $x,y\in\matN$.

Conversely suppose $d=12(x^2+3xy+3y^2)$ for $x,y\in\matN$, set $n=2x$ and $m=2(x+3y)$
and define $\ftil(z)=((n+m)-\omega m)\myprod z+\omega$. Then
$\ftil\left(\Btil^{(3)}_{0,0}\right)=\Ctil^{(6)}_{1,0}$ and
$\ftil\left(\Atil^{(2)}_{1,0}\right)$ is some $\Ctil^{(6)}_*$ by the above calculations, while
$$\ftil\left(\Ctil^{(6)}_{0,1}\right)=\ftil(\omega)=(2x+3y)+(2x+1)\omega$$
which is some $\Ctil^{(6)}_*$ because $(2x+3y)-(2x+1)\equiv 2\ ({\rm mod}\ 3)$,
hence $\ftil$ induces a realization of the candidate.
\end{proof}

\begin{thm}[case (3) in Fig.~\ref{gen:eucl:fig}]\label{333:to:333:thm}
The candidate surface branched covers inducing $S(3,3,3)\dotsto^{d:1} S(3,3,3)$ are
$$S\argdotstoqui{3k+1:1}{(3,\dots,3,1),(3,\dots,3,1),(3,\dots,3,1)}S\qquad
S\argdotstoqui{3k+3:1}{(3,\dots,3),(3,\dots,3),(3,\dots,3,1,1,1)}S$$
for $k\geqslant 1$, and they are realizable if and only if, respectively
\begin{itemize}
\item $d=x^2+xy+y^2$ with $x,y\in\matN$ and $x\nequiv y\ ({\rm mod}\ 3)$;
\item $d=3(x^2+3xy+3y^2)$ with $x,y\in\matN$.
\end{itemize}
\end{thm}

\begin{proof}
The possible covering instructions are
$$3\dotsto 3,\ 3\dotsto 3,\ 3\dotsto 3,\qquad
(3,3) \dotsto 3,\ 3\dotsto 3,\qquad
(3,3,3)\dotsto 3.$$
The second one is not induced by any candidate surface branched cover, and
the other two are induced by the candidates in the statement.

We follow again the same scheme, using the notation
of~(\ref{333:lifts:eq}). If $\ftil(z)=\lambda\myprod z+\mu$
realizes a candidate then, as in the previous proof,
Lemma~\ref{lattice:to:lattice:lem} implies that $\lambda=(n+m)-\omega m$
for $n,m\in\matZ$, and $d=n^2+nm+m^2$. Moreover from the symmetry of
$S(3,3,3)$ we can assume $f\left(A^{(3)}\right)=A^{(3)}$, whence
$\ftil\left(\Atil^{(3)}_{0,0}\right)=\Atil^{(3)}_{0,0}$, namely $\mu=0$.

For the first candidate we have $f\left(B^{(3)}\right)=B^{(3)}$ up
to symmetry, whence $\ftil\left(\Btil^{(3)}_{1,0}\right)=\lambda$ is some
$\Btil^{(3)}_{a,b}$ with $a-b\equiv 1\ ({\rm mod}\ 3)$,
therefore $n-m\equiv 1\ ({\rm mod}\ 3)$. Exactly as in the previous proof
we deduce that $d=x^2+xy+y^2$ with $x,y\in\matN$ and $x\nequiv y\ ({\rm mod}\ 3)$.
The converse is proved as above. Switching signs if necessary we assume
$x-y\equiv 1\ ({\rm mod}\ 3)$, we set $\ftil(z)=((x+y)-\omega y)\cdot z$
and note that $\ftil\left(\Ctil^{(3)}_{0,1}\right)=\ftil(\omega)=y+\omega x$
is $\Ctil^{(6)}_{y,x}$ because $y-x\equiv 2\ ({\rm mod}\ 3)$.

For the second candidate $f\left(B^{(3)}\right)=A^{(3)}$,
whence $\ftil\left(\Btil^{(3)_{1,0}}\right)=\lambda$ is some
$\Atil^{(3)}_*$, which implies that $n\equiv m\ ({\rm mod}\ 3)$.
Setting $x=n$ and $y=(m-n)/3$ we see that $x,y\in\matZ$ and $d=3(x^2+3xy+y^2)$,
and the conclusion is as usual.
\end{proof}

\begin{thm}[case (4) in Fig.~\ref{gen:eucl:fig}]\label{2222:to:2222:thm}
The candidate surface branched covers inducing some $S(2,2,2,2)\dotsto^{d:1} S(2,2,2,2)$ are
$$\begin{array}{l}
S\argdotstoset{2k+1:1}{(2,\dots,2,1),(2,\dots,2,1),(2,\dots,2,1),(2,\dots,2,1)}S\\
S\argdotstoset{2k+2:1}{(2,\dots,2,1,1),(2,\dots,2,1,1),(2,\dots,2),(2,\dots,2)}S\\
S\argdotstoset{2k+4:1}{(2,\dots,2,1,1,1,1),(2,\dots,2),(2,\dots,2),(2,\dots,2)}S
\end{array}$$
for $k\geqslant 1$. The first two are always realizable, the last one is if
and only if $d$ is a multiple of $4$.
\end{thm}

\begin{proof}
The first assertion is easy (but note that now the Riemann-Hurwitz formula
cannot be used in its simplified form~(\ref{RH:SS3:eq}), it reads
$\ell(\Pi)=2d+2$). The second assertion is proved as usual, except that
we have to deal with the flexibility of $S(2,2,2,2)$. We assume that a map
$f:S(2,2,2,2)\myto^{d:1} S(2,2,2,2)$
realizing some candidate exists and we put on the target $S(2,2,2,2)$
the structure $\pi$ defined by $\Gamma^{1,1}_{(2,2,2,2)}$.
Then we deduce from Lemma~\ref{Eucl:X:then:Eucl:Xtil:prop}
that there exists a structure $\pitil$ on the
source $S(2,2,2,2)$ also with area $2$,
and $\ftil:\matE\to\matE$, such that
$\ftil(z)=\lambda\myprod z+\mu$ with $\pi\compo\ftil=f\compo\pitil$ and
$d=|\lambda|^2$. Then $\pitil$ is defined by some $\Gamma^{s,t}_{(2,2,2,2)}$.

We first note that by the symmetry of $S(2,2,2,2)$ we can assume $\mu=0$.
Then we apply Lemma~\ref{lattice:to:lattice:lem}.
Since
$$\Lambda^{1,1}_{(2,2,2,2)}=\langle\tau_2,\tau_{2i}\rangle,\qquad
\left\langle\tau_{2is},\tau_{2\left(\frac{1}{s}+it\right)}\right\rangle$$
there exist $m,n,p,q \in \matZ$ such that
\begin{equation}\label{2222:to:2222:eq}
\left\{\begin{array}{l}
\lambda\myprod is=n+im\\
\lambda\myprod \left(\frac{1}{s}+it\right)=p+iq
\end{array}
\right.
\end{equation}
and some easy computations show that all
the relevant quantities can be determined
explicitly in terms of $n,m,p,q$, namely:
$$s=\sqrt{\frac{n^2+m^2}{pm-qn}},\qquad
t=\frac{sp}n-\frac m{ns},\qquad
\lambda=\frac 1s(m-in)$$
so in particular $d=|\lambda|^2= pm-qn$.
Note also that equations~(\ref{2222:to:2222:eq})
already give us also the images of the lifts of the cone points.

For the first candidate we have $d=4a\pm 1$ for some $a\geqslant 1$, we set
$n=2a,\ m=1,\ p=\pm1,\ q=-2$, we compute $s,t,\lambda$ as above and we see that the
corresponding map $\ftil$ induces a realization of the candidate, because
$$\begin{array}{ll}
\ftil\left(\Atil^{(2)}\right)=\Atil^{(2)}_{0,0} &
\ftil\left(\Btil^{(2)}\right)=\Btil^{(2)}_{\pm1,-2} \\
\ftil\left(\Ctil^{(2)}\right)=\Ctil^{(2)}_{2a\pm 1,-1} &
\ftil\left(\Dtil^{(2)}\right)=\Dtil^{(2)}_{2a,1}.
\end{array}$$

For the second candidate we have $d=4a+1\pm1$ for some $a\geqslant 1$, we set
$n=2a,\ m=1\pm 1,\ p=1,\ q=-2$
we compute $s,t,\lambda$ as above and we see that the
corresponding map $\ftil$ induces a realization of the candidate, because
$$\begin{array}{ll}
\ftil\left(\Atil^{(2)}\right)=\Atil^{(2)}_{0,0} &
\ftil\left(\Btil^{(2)}\right)=\Btil^{(2)}_{1,-2} \\
\ftil\left(\Ctil^{(2)}\right)=\Btil^{(2)}_{2a+1,\pm 1-1} &
\ftil\left(\Dtil^{(2)}\right)=\Atil^{(2)}_{2a,1\pm1}.
\end{array}$$

For the last candidate each lift of a cone point has some
$\Atil^{(2)}_*$ as its image, therefore $n,m,p,q$ must all be even,
which implies that $d$ is a multiple of $4$, as prescribed in the statement.
Conversely if $d=4a$ for $a>1$ we set
$n=m=q=2,\ p=2(a+1)$
we compute $s,t,\lambda$ as above and we see that the
corresponding map $\ftil$ induces a realization of the candidate, because
$$\begin{array}{ll}
\ftil\left(\Atil^{(2)}\right)=\Atil^{(2)}_{0,0} &
\ftil\left(\Btil^{(2)}\right)=\Atil^{(2)}_{2(a+1),2} \\
\ftil\left(\Ctil^{(2)}\right)=\Atil^{(2)}_{2(a+2),4} &
\ftil\left(\Dtil^{(2)}\right)=\Atil^{(2)}_{2,2}.
\end{array}$$
The proof is complete.
\end{proof}

\begin{thm}[case (5) in Fig.~\ref{gen:eucl:fig}]\label{333:to:236:thm}
The candidate surface branched covers inducing some $S(3,3,3)\dotsto^{d:1} S(2,3,6)$ are
$$\begin{array}{ll}
S\argdotstoqui{6k:1}{(2,\dots,2),(3,\dots,3,1,1,1),(6,\dots,6)}S\quad &
S\argdotstoqui{6k+2:1}{(2,\dots,2),(3,\dots,3,1,1),(6,\dots,6,2)}S \\
S\argdotstoqui{6k+4:1}{(2,\dots,2),(3,\dots,3,1),(6,\dots,6,2,2)}S\quad &
S\argdotstoqui{6k+6:1}{(2,\dots,2),(3,\dots,3),(6,\dots,6,2,2,2)}S
\end{array}$$
for $k\geqslant 1$, and they are realizable, respectively:
\begin{itemize}
\item if and only if $d=6(x^2+3xy+3y^2)$ for $x,y\in\matN$;
\item if and only if $d=2(x^2+xy+y^2)$ for $x,y\in\matN$ and
$x\nequiv y\ ({\rm mod}\ 3)$;
\item never;
\item if and only if $d=6(x^2+3xy+3y^2)$ for $x,y\in\matN$.
\end {itemize}
\end{thm}

\begin{proof}
The first assertion is easy. For the second one we proceed as above, except that now
the Euclidean structure $\pitil$ on $S(3,3,3)$ is not that we have fixed above, because
its area should be $\frac{\sqrt{3}}4$ rather than $\frac{\sqrt{3}}2$, so the triangle
$\Delta(3,3,3)$ must be rescaled by a factor $1/\sqrt{2}$. The lattices to which
we can apply Lemma~\ref{lattice:to:lattice:lem} are therefore
$$\frac 1{\sqrt2}\cdot\Lambda_{(3,3,3)}=\left\langle \tau_{i\sqrt{\frac32}},\tau_{\frac{3+i\sqrt3}{2\sqrt2}}\right\rangle,\qquad
\Lambda_{(2,3,6)}=\left\langle \tau_{i\sqrt3},\tau_{\frac{3+i\sqrt3}2}\right\rangle.$$
As in the proof of Theorem~\ref{236:to:236:thm} (except for the new factor) we deduce that
$$\lambda=\sqrt2\myprod((n+m)-\omega m),\qquad
d=|\lambda|^2=2(n^2+nm+m^2).$$
Therefore $\ftil$ maps the lifts of the cone points of $S(3,3,3)$ to
$$\ftil(0)=\mu,\quad
\ftil\left(\frac 1{\sqrt2}\right)=(n+m)-\omega m+\mu,\quad
\ftil\left(\frac {\omega}{\sqrt2}\right)=m+n\omega+\mu.$$
For the first candidate all these points should be some $\Btil^{(3)}_*$
from~(\ref{236:lifts:eq}), so we can assume $\mu=0$ and
$$n+m-(-m)\nequiv 2\ ({\rm mod}\ 3),\quad m-n \nequiv 2\ ({\rm mod}\ 3)\quad
\Rightarrow\qquad n\equiv m\ ({\rm mod}\ 3).$$
Setting $x=n$ and $y=(m-n)/3$ we then see that $d=6(x^2+3xy+3y^2)$ for $x,y\in\matZ$,
and as above we can reduce to $x,y\in\matN$, so $d$ has the appropriate form.
The converse follows from the same computations: if
$d=6(x^2+3xy+3y^2)$ we set $n=x$ and $m=x+3y$ and we see that the corresponding
$\ftil$ realizes the candidate.

For
the second candidate again $\mu=0$ and, by the symmetry of $S(3,3,3)$,
we can assume $1/\sqrt2$ is mapped to some $C^{(6)}_*$, namely
$n-m\equiv 2\ ({\rm mod}\ 3)$, so in particular $n\nequiv m\ ({\rm mod}\ 3)$.
Therefore $d=2(x^2+xy+y^2)$ for some $x,y\in\matZ$ with
$x\nequiv y\ ({\rm mod}\ 3)$, and once again restricting to $x,y\in\matN$ makes
no difference, so $d$ has the prescribed form. The construction is easily
reversible because if $n-m\equiv 2\ ({\rm mod}\ 3)$ then $m-n \nequiv 2\ ({\rm mod}\ 3)$,
which also proves that the third candidate is never realizable.

For the last candidate we can assume $\mu=\omega$, and
$$
\begin{array}{l}
(n+m)-(1-m)\equiv 2\ ({\rm mod}\ 3)\\
m-(n+1) \equiv 2\ ({\rm mod}\ 3)
\end{array}
\quad
\Rightarrow\qquad n\equiv m\ ({\rm mod}\ 3)$$
and we conclude as for the first candidate.
\end{proof}

\begin{thm}[case (6) in Fig.~\ref{gen:eucl:fig}]\label{2222:to:244:thm}
The candidate surface branched covers inducing some $S(2,2,2,2)\dotsto^{d:1} S(2,4,4)$
are
$$\begin{array}{ll}
S\argdotstosex{4k+4:1}{(2,\dots,2,1,1,1,1),(4,\dots,4),(4,\dots,4)}S\quad &
S\argdotstosex{4k+4:1}{(2,\dots,2,1,1),(4,\dots,4,2,2),(4,\dots,4)}S\\
S\argdotstosex{4k+2:1}{(2,\dots,2,1,1),(4,\dots,4,2),(4,\dots,4,2)}S\quad &
S\argdotstosex{4k+4:1}{(2,\dots,2),(4,\dots,4,2,2),(4,\dots,4,2,2)}S\\
S\argdotstosex{4k+6:1}{(2,\dots,2),(4,\dots,4,2,2,2),(4,\dots,4,2)}S\quad &
S\argdotstosex{4k+8:1}{(2,\dots,2),(4,\dots,4,2,2,2,2),(4,\dots,4)}S
\end{array}$$
for $k\geqslant 1$. The first four are always realizable, the fifth
one is never, and the last one is if and only if $d$ is a multiple of $8$.
\end{thm}

\begin{proof}
Again we leave the first assertion to the reader and we proceed with the
customary scheme. Since the area of the structure we have chosen on $S(2,4,4)$ is $1$,
on $S(2,2,2,2)$ we will have a structure generated by the rotations of angle $\pi$
around points
$$0,\quad \frac1{2s}+it,\quad \frac1{2s}+i(s+t),\quad is$$
with $s,t\in\matR$ and $s>0$. The lattices to which we must apply
Lemma~\ref{lattice:to:lattice:lem} are therefore
$\left\langle\tau_{2is},\tau_{\frac{1}{s}+2it}\right\rangle$ and
$\Lambda_{(2,4,4)}=\left\langle\tau_2,\tau_{2i}\right\rangle$, so
\begin{displaymath}
\left\{\begin{array}{l}
\lambda\myprod 2is= 2(n+im)\\
\lambda\myprod (\frac{1}{s}+2it)=2(p+iq)
\end{array}
\right.
\end{displaymath}
for some $n,m,p,q\in\matZ$. Whence, after easy computations,
$$s=\sqrt{\frac{n^2+m^2}{2(pm-qn)}},\qquad t=\frac{sp}n-\frac m{2sn},\qquad
\lambda=\frac{m-in}s.$$
In particular $d=2(pm-qn)$ and the images of the lifts of the cone points of $S(2,2,2,2)$ are
$$\begin{array}{ll}
\ftil(0)=\mu\quad &
\ftil\left(\frac1{2s}+it\right)=p+iq+\mu\\
\ftil(is)=n+im+\mu\quad &
\ftil\left(\frac1{2s}+i(s+t)\right)=(p+n)+i(q+m)+\mu.
\end{array}$$

The first four candidates are realized respectively with the following choices of $n,m,p,q,\mu$:
$$\begin{array}{c|c|c|c|c}
n & m & p & q & \mu \\ \hline
k+1   & k+1    & 1      & -1      & 1  \\
k       & k+1     & 2     & 0 & 0 \\
k       & k+1      & 1      & -1    & 0 \\
k       & k+1      & 2    & 0    & 0
\end{array}$$

The fifth candidate is always exceptional because we can suppose $\mu=0$ and
hence we should have that two of the pairs
$$(p,q),\quad (n,m),\quad (p+n,q+m)$$
consist of even numbers and the third one consists of odd numbers, which is impossible.

For the last candidate we have that $p,q,n,m$ must all be even, so $d=2(pm-qn)$ is a multiple of $8$.
Conversely if $d=8h$ we can realize the candidate with $n=q=0$, $m=2$ and $p=2h$.
\end{proof}

\begin{thm}[case (7) in Fig.~\ref{gen:eucl:fig}]\label{2222:to:236:thm}
The candidate surface branched covers inducing $S(2,2,2,2)\dotsto^{d:1} S(2,3,6)$ are
$$\begin{array}{ll}
S\argdotstosex{6k:1}{(2,\dots,2,1,1,1,1),(3,\dots,3),(6,\dots,6)}S\quad &
S\argdotstosex{6k+3:1}{(2,\dots,2,1,1,1),(3,\dots,3),(6,\dots,6,3)}S\\
S\argdotstosex{6k+6:1}{(2,\dots,2,1,1),(3,\dots,3),(6,\dots,6,3,3)}S\quad &
S\argdotstosex{6k+9:1}{(2,\dots,2,1),(3,\dots,3),(6,\dots,6,3,3,3)}S\\
S\argdotstosex{6k+12:1}{(2,\dots,2),(3,\dots,3),(6,\dots,6,3,3,3,3)}S\quad &
\end{array}$$
for $k\geqslant 1$. The first three are always realizable, the fourth
one is never, and the last one is if and only if $d$ is a multiple of $12$.
\end{thm}

\begin{proof}
Once again we leave the first assertion to the reader and we follow the usual
scheme. Since the area of $S(2,3,6)$ is $\sqrt{3}/4$,
on $S(2,2,2,2)$ we will have a structure generated by the rotations of angle $\pi$
around points
$$0,\quad \frac{\sqrt3}{8s}+it,\quad \frac{\sqrt3}{8s}+i(s+t),\quad is$$
and we apply Lemma~\ref{lattice:to:lattice:lem} to
$\left\langle\tau_{2is},\tau_{\frac{\sqrt3}{8s}+2it}\right\rangle$ and
$\Lambda_{(2,3,6)}=
\left\langle \tau_{i\sqrt3},\tau_{\frac{3+i\sqrt3}2}\right\rangle$, so for some $n,m,p,q\in\matZ$ we have
\begin{displaymath}
\left\{\begin{array}{l}
\lambda\myprod 2is=ni\sqrt3+m\frac{3+i\sqrt3}{2}\\
\lambda\myprod \left(\frac{\sqrt3}{4s}+2it\right)=pi\sqrt3+q\frac{3+i\sqrt3}{2}
\end{array}
\right.
\end{displaymath}
whence, after some calculations
\begin{eqnarray*}
& & s=\frac 12\sqrt{\frac{n^2+nm+m^2}{qn-pm}},\qquad t=\frac{qs}m-\frac{m+2n}{8ms},\qquad \\
& & \lambda=\frac{\sqrt3(m+2n)-3im}{4s}=\frac{\sqrt3}{2s}\myprod((n+m)-m\omega)
\end{eqnarray*}
so in particular $d=|\lambda|^2=3(qn-pm)$. Moreover the following relations
will readily allow us to determine the images under $\ftil$
of the lifts of the cone points of $S(2,2,2,2)$:
\begin{eqnarray*}
& & \lambda\myprod is=\frac 12((m-n)+(m+2n)\omega),\\
& & \lambda\myprod \left(\frac{\sqrt3}{8s}+it\right)=
\frac 12((q-p)+(q+2p)\omega).
\end{eqnarray*}

For the first candidate we choose $\mu=\frac 12$, $p=q=2$, $n=k+1$ and $m=1$.
The corresponding $\ftil$ induces a realization because
$d=6k=3(qn-pm)$ and
the images of the
cone points are
$$\begin{array}{ll}
\frac 12 (1+0\omega),\phantom{\Big|}\ & \frac 12((m-n+1)+(m+2n)\omega),\\
\frac 12 ((q-p+1)+(q+2p)\omega),\ & \frac 12((m+q-n-p+1)+(m+q+2n+2p)\omega)
\end{array}$$
which are easily recognized to all have the form $\frac 12(a+b\omega)$ with $a,b$ not both even and
$a-b\equiv 1\ ({\rm mod}\ 3)$, so they equal some $\Atil^{(2)}_*$.

For the second candidate we choose $\mu=\frac 12$, $n=2$, $m=1$ and
$$q=k,\ p=-1\quad {\rm if}\ k\equiv 1\ ({\rm mod}\ 2),\qquad
q=k+1,\ p=1\quad {\rm if}\ k\equiv 0\ ({\rm mod}\ 2).$$
Then $d=6k+3=3(qn-pm)$ and the images of the cone points are as before,
but now the first three are some $\Atil^{(2)}_*$, while the last one
has the form $a+b\omega$ with $a,b\in\matZ$ and
$a-b\equiv 2\ ({\rm mod}\ 3)$, so it is some $\Ctil^{(6)}_*$,
so $\ftil$ induces a realization of the candidate.

For the third candidate we choose $\mu=\frac 12$, $m=q=2$ and
$$n=k,\ p=-1\quad {\rm if}\ k\equiv 0\ ({\rm mod}\ 2),\qquad
n=k+1,\ p=1\quad {\rm if}\ k\equiv 1\ ({\rm mod}\ 2).$$
Then $d=6k+6=3(qn-pm)$ and now the
first two images are some $\Atil^{(2)}_*$ and the last two are some $\Ctil^{(6)}_*$,
so $\ftil$ induces a realization of the candidate.

For the fourth candidate we can once again suppose $\mu=\frac 12$. Since the images of the
last three cone points must be some $\Ctil^{(6)}_*$ we deduce that
$m,q,m+q$ should be even and hence $n,p,n+p$ should be odd,
which is impossible.

Turning to the last candidate, we can suppose $\mu=\omega$. Then the images of the cone points are
$$\begin{array}{ll}
1+0\omega,\phantom{\Big|}\ & \frac 12((m-n)+(m+2n+2)\omega),\\
\frac 12 ((q-p)+(q+2p+2)\omega),\ & \frac 12((m+q-n-p)+(m+q+2n+2p+2)\omega)
\end{array}$$
and they must all be some $\Ctil^{(6)}_*$, so $n,m,p,q$ should all be even.
Therefore $d=3(qn-pm)$ is a multiple of $12$. Conversely, if $d=12h+12$ we realize the candidate
with the choice $q=2h$, $n=m=2$ and $p=-2$.
\end{proof}

\paragraph{Congruences and density}
We close this section with some elementary remarks about
the numbers appearing in our main results, thus explaining
in full detail why Theorems~\ref{244:main:thm}
to~\ref{333:main:thm} are implied by Theorems~\ref{244:to:244:thm}
to~\ref{333:to:333:thm}. First of all we have:
$$\begin{array}{l}
\big\{d\in\matN:\ d=x^2+y^2\ {\rm for}\ x,y\in\matN,\ x\nequiv y\ ({\rm mod}\ 2)\big\} \\
= \big\{d\in\matN:\ d\equiv 1\ ({\rm mod}\ 4),\ d=x^2+y^2\ {\rm for}\ x,y\in\matN\big\},\phantom{\Big|}\\
\big\{d\in\matN:\ d=x^2+xy+y^2\ {\rm for}\ x,y\in\matN\ {\rm not\ both\ even},\ x\nequiv y\ ({\rm mod}\ 3)\big\} \\
= \big\{d\in\matN:\ d\equiv 1\ ({\rm mod}\ 6),\ d=x^2+xy+y^2\ {\rm for}\ x,y\in\matN\big\},\phantom{\Big|}\\
\big\{d\in\matN:\ d=x^2+xy+y^2\ {\rm for}\ x,y\in\matN,\ x\nequiv y\ ({\rm mod}\ 3)\big\} \\
= \big\{d\in\matN:\ d\equiv 1\ ({\rm mod}\ 3),\ d=x^2+xy+y^2\ {\rm for}\ x,y\in\matN\big\}.\phantom{\Big|}
\end{array}$$
Moreover the statement made in the Introduction that in Theorems~\ref{244:main:thm}
to~\ref{333:main:thm} the
realizable degrees have zero asymptotic density means the following:
$$\begin{array}{l}
\lim\limits_{n\to\infty}\frac 1n \myprod \# \big\{d\in \matN:\ d\leqslant n,\ d=x^2+y^2\ {\rm for}\ x,y\in\matN\big\} = 0, \\
\lim\limits_{n\to\infty}\frac 1n \myprod \# \big\{d\in \matN:\ d\leqslant n,\ d=x^2+xy+y^2\ {\rm for}\ x,y\in\matN\big\} = 0.
\end{array}$$

\section{Hyperbolic triangular orbifolds}\label{neg:chi:sec}
There is one crucial geometric fact underlying the proofs of our
main results for the case of positive and zero Euler
characteristic. Namely, in these cases the geometry (if any) of an
orbifold with cone points is \emph{rigid} (up to rescaling), with the single
exception of $S(2,2,2,2)$, where the space of moduli is easy
to compute anyway. Turning to the case of negative Euler
characteristic, one knows that a hyperbolic 2-orbifold
is rigid if and only if it is \emph{triangular}, namely if it is
based on the sphere and it has precisely three cone points. In this
section we will show that only very few candidate surface branched
covers induce candidate covers between hyperbolic triangular
2-orbifolds:

\begin{thm}\label{hyp:tria:list:thm}
The candidate surface branched covers
inducing candidate covers between triangular hyperbolic $2$-orbifolds are precisely:
$$\begin{array}{lll}
S\argdotstoqui{6:1}{(5,1),(4,1,1),(2,2,2)}S\quad &
S\argdotstoqui{8:1}{(5,1,1,1),(4,4),(2,\ldots,2)}S\quad &
S\argdotstoqui{8:1}{(7,1),(3,3,1,1),(2,\ldots,2)}S\quad \\
S\argdotstoqui{9:1}{(7,1,1),(3,3,3),(2,\ldots,2,1)}S\quad &
S\argdotstoqui{10:1}{(8,1,1),(3,3,3,1),(2,\ldots,2)}S\quad &
S\argdotstoqui{12:1}{(8,2,1,1),(3,\ldots,3),(2,\ldots,2)}S\quad \\
S\argdotstoqui{12:1}{(9,1,1,1),(3,\ldots,3),(2,\ldots,2)}S\quad &
S\argdotstoqui{16:1}{(7,7,1,1),(3,\ldots,3,1),(2,\ldots,2)}S\quad &
S\argdotstoqui{24:1}{(7,7,7,1,1,1),(3,\ldots,3),(2,\ldots,2)}S.\quad
\end{array}$$
\end{thm}

This result means that the geometric techniques employed above would
require a substantial extension to be relevant for the hyperbolic case.
Leaving this for the future, we will show here that for
the candidate covers of Theorem~\ref{hyp:tria:list:thm}
the geometric approach is not even necessary, since realizability
can be fully analyzed using a completely different technique, namely Grothendieck's
\emph{dessins d'enfant}~\cite{Groth,Wolfart}, already exploited in~\cite{PePe1}.
We will show the following:

\begin{prop}\label{hyp:tria:real:prop}
Among the candidate covers of Theorem~\ref{hyp:tria:list:thm}, the second and the eighth are
exceptional and all other ones are realizable.
\end{prop}

Let us now establish the results we have stated.
The first proof requires the analysis
of quite a few cases, some of which we will leave to the reader.

\dimo{hyp:tria:list:thm}
Our argument is organized in three steps:
\begin{itemize}
\item[(I)] Analysis of the relevant surface candidate covers with degree $d\leqslant 11$;
\item[(II)] Restrictions on the base of the induced candidate cover for $d\geqslant 12$;
\item[(III)] More restrictions on the cover and conclusion for $d\geqslant 12$.
\end{itemize}

\bigskip

\noindent\textsc{Step I.}
If $\Pi$ is a partition of an integer $d$, let us denote by $\ell(\Pi)$ its length
(as above), and by $c(\Pi)$ the number of entries in
$\Pi$ which are different from ${\rm l.c.m.}(\Pi)$.
To induce a candidate cover between triangular 2-orbifolds (regardless of the geometry),
a candidate surface branched cover of degree $d\geqslant 2$ must have the following properties:
\begin{itemize}
\item The number of branching points is 3;
\item If the partitions of $d$ are $\Pi_1,\Pi_2,\Pi_3$ then $c(\Pi_1)+c(\Pi_2)+c(\Pi_3)=3$.
\end{itemize}
To list all such candidate covers for a given $d$ then one has to:
\begin{itemize}
\item List all the partitions $\Pi$ of $d$ with $c(\Pi)\leqslant 3$;
\item Find all possible triples $(\Pi_1,\Pi_2,\Pi_3)$ of partitions with
$\ell(\Pi_1)+\ell(\Pi_2)+\ell(\Pi_3)=d+2$ and $c(\Pi_1)+c(\Pi_2)+c(\Pi_3)=3$.
\end{itemize}
We have done this for $2\leqslant d\leqslant 11$ and then we have singled out the
candidate covers inducing hyperbolic 2-orbifold covers, getting
the first five items of the statement.
To illustrate how this works we will spell out here only
the case $d=8$. The partitions $\Pi$ of $8$ with $c(\Pi)\leqslant 3$ are those described in
Table~\ref{c3:part:8:tab}, with the corresponding values of $\ell$ and $c$.

\begin{table}
\begin{center}
\begin{tabular}{c||c|c|c|c|c}
$\Pi$       & (8)       & (6,1,1)   & (5,1,1,1) & (4,2,2)   & (3,3,1,1)       \\ \hline
$\ell$      & 1         & 3         & 4         &  3        &  4        \\ \hline
$c$         & 0         & 2         & 3         &  2        &  2        \\ \hline\hline
$\Pi$       & (7,1)     & (5,3)     & (4,4)     & (4,2,1,1) & (2,2,2,2)       \\ \hline
$\ell$      &  2        & 2         & 2         &  4        &  4        \\ \hline
$c$         &  1        & 2         & 0         &  3        &  0        \\ \hline\hline
$\Pi$       & (6,2)     & (5,2,1)   & (4,3,1)   & (3,3,2)   & (2,2,2,1,1)       \\ \hline
$\ell$      &  2        & 3         & 3         &  3        &  5        \\ \hline
$c$         &  1        & 3         & 3         &  3        &  2        \\
\end{tabular}
\end{center}
\mycap{The partitions $\Pi$ of $8$ with $c(\Pi)\leqslant 3$\label{c3:part:8:tab}}
\end{table}
The triples of such partitions such that $\ell$ and $c$ sum up to $10$ and $3$ respectively
are shown in Table~\ref{tria:8:tab}, together with the associated candidate orbifold cover and
its geometric type.
\begin{table}
\begin{center}
\begin{tabular}{c|c|c|c|c}
\multicolumn{3}{c|}{$\Pi_1,\Pi_2,\Pi_3$} & Associated cover & Geometry \\ \hline\hline
(4,2,1,1) & (4,4) & (2,2,2,2) & $S(2,4,4)\dotsto S(2,4,4)$ & $\matE$ \\ \hline
(5,1,1,1) & (4,4) & (2,2,2,2) & $S(5,5,5)\dotsto S(2,4,5)$ & $\matH$ \\ \hline
(6,2) & (3,3,1,1) & (2,2,2,2) & $S(3,3,3)\dotsto S(2,3,6)$ & $\matE$ \\ \hline
(7,1) & (3,3,1,1) & (2,2,2,2) & $S(3,3,7)\dotsto S(2,3,7)$ & $\matH$ \\
\end{tabular}
\end{center}
\mycap{Triples of partitions of $8$ inducing candidate covers between
triangular orbifolds\label{tria:8:tab}}
\end{table}
So we get the second and third item in the statement.

\bigskip

\noindent\textsc{Step II.}
Let us denote by $\Xtil\dotsto^{d:1} X$ a candidate orbifold cover with $d\geqslant 12$
and hyperbolic $\Xtil=S(\alpha,\beta,\gamma)$ and $X=S(p,q,r)$. Since
$$0<-\chiorb(\Xtil)=1-\left(\frac1\alpha+\frac1\beta+\frac1\gamma\right)<1$$
and $\chiorb(\Xtil)=d\myprod\chiorb X$, we deduce that
$$0<-\chiorb(X)=1-\left(\frac1p+\frac1q+\frac1r\right)<\frac1{12}
\quad\Rightarrow\quad
\frac{11}{12}<\left(\frac1p+\frac1q+\frac1r\right)<1.$$
Assuming $p\leqslant q\leqslant r$ it is now very easy to check that
the last inequality is satisfied only for
$p=2$, $q=3$, $7\leqslant r\leqslant 11$ and for $p=2$, $q=4$, $r=5$.

\bigskip

\noindent\textsc{Step III.}
If $\Xtil\dotsto X$ is a candidate 2-orbifold cover with hyperbolic
$\Xtil=S(\alpha,\beta,\gamma)$ and $X=S(p,q,r)$
then the following must happen:
\begin{itemize}
\item[(a)] Each of $\alpha,\beta,\gamma$ must be a divisor of some
element of $\{p,q,r\}$;
\item[(b)] $\frac{\chiorb(\Xtil)}{\chiorb(X)}$ must be an integer $d$;
\item[(c)] There must exist three partitions of $d$ inducing $\Xtil\dotsto X$.
\end{itemize}
Successively imposing these conditions with each of the 5 orbifolds $X$
coming from Step II and restricting to $d\geqslant 12$ we have found the
last four items in the statement. Again we only
spell out here one example, leaving the other ones to the reader.
Let $X$ be $S(2,3,8)$. Then the relevant hyperbolic $\Xtil$'s
according to (a), excluding $X$ itself, are
$$\begin{array}{ccccc}
S(2,4,8) & S(3,3,4) & S(2,8,8) & S(3,3,8) & S(3,8,8) \\
S(3,4,4) & S(4,4,4) & S(4,4,8) & S(4,8,8) & S(8,8,8)
\end{array}$$
and $d=\frac{\chiorb(\Xtil)}{\chiorb(X)}$ is always integer in this case,
so point (b) is not an issue.
However $d\leqslant 11$ in all cases but the last two (for instance,
the case $\Xtil=S(3,8,8)$ corresponds to the fifth item in the statement).
For $\Xtil=S(4,8,8)$ we have $d=12$
and taking the partitions of $12$ as the sixth item in the statement
we see that the induced orbifold cover is indeed
$S(4,8,8)\dotsto S(2,3,8)$. For $\Xtil=S(8,8,8)$ we have
$d=15$ and it is impossible to find partitions of $15$ inducing the candidate cover,
because the cone point of order 2 in $X$, being covered by ordinary points of $\Xtil$
only, should require a partition consisting of $2$'s only, which
cannot exist because $15$ is odd.

Carrying out the same analysis one gets the last two items in the statement
for $X=S(2,3,7)$, the seventh item
for $X=S(2,3,9)$, and nothing new for the other $X$'s.
This concludes Step III and the proof.
\finedimo

As already announced, the next argument is based on a technique different from
those used in the rest of this paper, namely
Grothendieck's \emph{dessins d'enfant}. We will not review this
tool here, addressing the reader to~\cite{PePe1}.

\begin{figure}
\centering
\includegraphics[width=11cm]{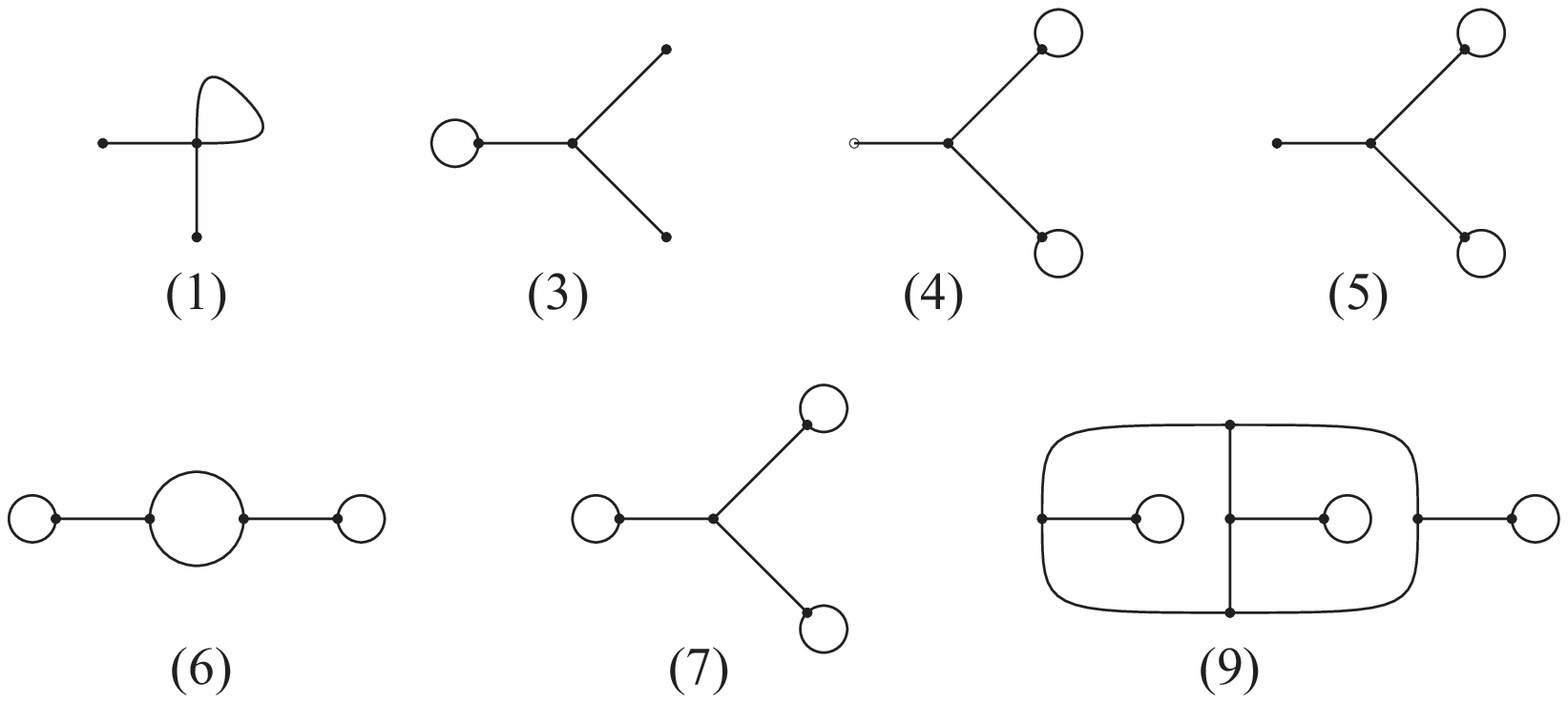}
\mycap{Dessins d'enfant for all candidate surfaced branched covers in
Theorem~\ref{hyp:tria:list:thm} except the second and the eighth\label{hyp_tria_OK_DDE}}
\end{figure}

\dimo{hyp:tria:real:prop}
Dessins d'enfant proving the realizability of all candidate covers
claimed to be realizable can be found in Fig.~\ref{hyp_tria_OK_DDE}.
The black vertices always correspond to the elements of the second partition in
Theorem~\ref{hyp:tria:list:thm}, and the white vertices to the entries of the
third partition, while the regions correspond to the
elements of the first partition. However 2-valent white vertices are never shown,
except for the single 1-valent one in case (4).

\begin{figure}
\centering
\includegraphics[width=3.5cm]{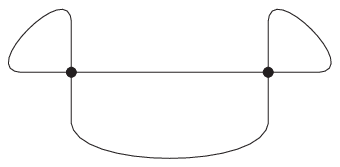}
\mycap{Exceptionality of $S\argdotstoqui{8:1}{(5,1,1,1),(4,4),(2,\ldots,2)}S$\label{hyp_tria_KO2_DDE}}
\end{figure}

Exceptionality of $S\argdotstoqua{8:1}{(5,1,1,1),(4,4),(2,\ldots,2)}S$
is easy: a dessin relative to partitions $(4,4)$ and $(2,2,2,2)$
with at least two outer regions of length 1 must be as shown in Fig.~\ref{hyp_tria_KO2_DDE},
so the third partition is $(4,2,1,1)$, not $(5,1,1,1)$.

\begin{figure}
\centering
\includegraphics[width=11cm]{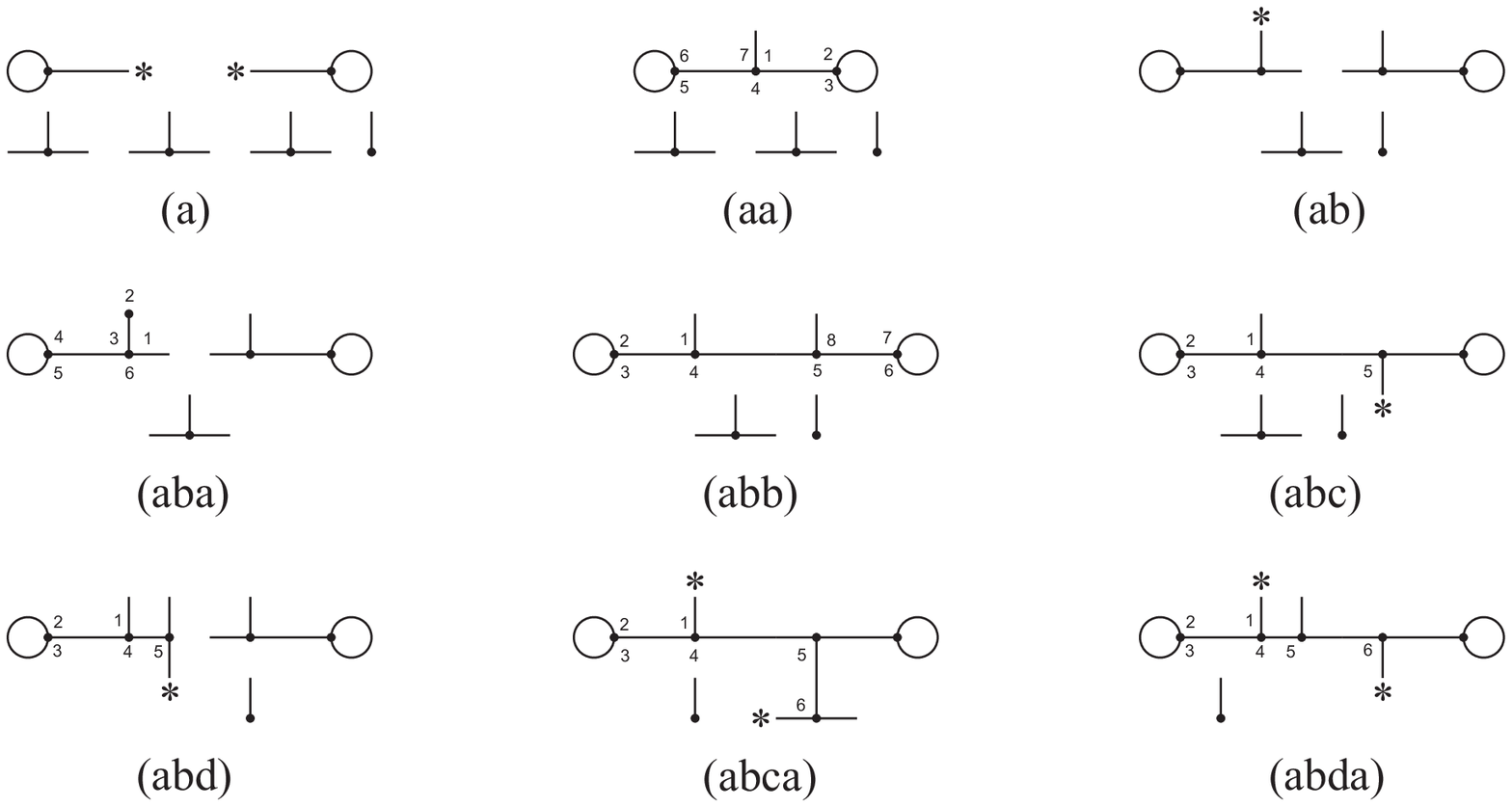}
\mycap{Exceptionality of $S\argdotstosex{16:1}{(7,7,1,1),(3,\ldots,3,1),(2,\ldots,2)}S$\label{hyp_tria_KO8_DDE}}
\end{figure}

For the exceptionality of $S\argdotstoqui{16:1}{(7,7,1,1),(3,\ldots,3,1),(2,\ldots,2)}S$
refer to Fig.~\ref{hyp_tria_KO8_DDE}.
Since it must contain two length-1 regions, a dessin realizing it
should be as in (a). The two marked germs of edges cannot be joined
together or to the 1-valent vertex, so they go either to the same 3-valent
vertex as in (aa) or to different 3-valent vertices as in (ab).
Case (aa) is impossible because there is a region with 7 vertices, which
will become more than 7 in the complete dessin. In case (ab) we examine
where the marked germ of edge could go, getting cases (aba) to (abd), always redrawn in
a more convenient way. Cases (aba) and (abb) are impossible because of long regions.
In cases (abc) and (abd) we examine where the marked edge could go in
order not to create regions of length 5 or longer than 7, and we see
that in both cases there is only one possibility, namely
(abca) and (abda). In both these cases, because of the region of length
already 6, the two marked germs of edges should go to one and the same
3-valent vertex, but there are no more available with two free germs
of edges, so again we cannot
complete the dessin in order to realize (8). Our argument is complete.
\finedimo

\vspace{.35cm}

\noindent
Dipartimento di Matematica\\
Universit\`a di Roma ``La Sapienza''\\
P.le Aldo Moro, 2\\
00185 ROMA, Italy\\
pascali@mat.uniroma1.it

\vspace{.35cm}

\noindent
Dipartimento di Matematica Applicata\\
Universit\`a di Pisa\\
Via Filippo Buonarroti, 1C\\
56127 PISA, Italy\\
petronio@dm.unipi.it

\end{document}